\documentclass{article}
\usepackage{amsmath, amsthm, amsfonts, amssymb}
\usepackage{hyperref}
\newtheorem{theorem}{Theorem}[section]
\newtheorem{lemma}[theorem]{Lemma}
\newtheorem{proposition}[theorem]{Proposition}
\newtheorem{corollary}[theorem]{Corollary}
\theoremstyle{definition}
\newtheorem{definition}[theorem]{Definition}
\newtheorem{example}[theorem]{Example}

\theoremstyle{remark}
\newtheorem{remark}[theorem]{Remark}
\numberwithin{equation}{section}

\renewcommand{\(}{\left(}
\renewcommand{\)}{\right)}

\renewcommand{\]}{\right]}

\begin{document}
\title{\bf Orlicz Spaces with Bicomplex Scalars}
\date{\text{ R. Kumar, K. Sharma, R. Tundup and S. Wazir}}
\vspace{0in}
\maketitle
$\textbf{Abstract.}$ 
In this paper we define bicomplex Orlicz space with hyperbolic valued Luxemburg norm and discussed some of their properties. We have also partially characterize an integral representation of a $\mathbb{D}$-valued convex function. Further we have shown that the spectrum of the unilateral shift operator on $l^2(\mathbb{BC})$ is the null cone $\mathcal{NC}$.\\\\
$\textbf{Keywords.}$ Modules, Hyperbolic Numbers, $\mathbb{D}$-metric, $\mathbb{D}$-Convex function, $\mathbb{D}$-Orlicz Function, $\mathbb{D}$-metric, $\mathbb{D}$-metric spaces, $\mathbb{D}$-Mulsielak Orlicz Function.\\\\
$\textbf{AMS  Subject Classification.}$ 30G35 47A10.

\begin{section} {Introduction}

\renewcommand{\thefootnote}{\fnsymbol{footnote}}
\footnotetext{The work of the third author is supported by Rajiv Gandhi National Fellowship(UGC), INDIA (No.F.17.1/2013-14/RGNF-ST-JAM-38257(SA-III)).}

\bigskip 
The study of hyperbolic numbers can be traced way back to 1848 from the work of James Cockle \cite{Co-1848} and then to 1893 by Lie and Scheffers \cite{Li-Sc93}. Hyperbolic number system has widely been studied for various reasons, one among which is its commutative property. Infact along with the set of complex numbers, the set of hyperbolic number were found to be the only real commutative Clifford algebra. The importance of hyperbolic numbers lies in the fact that the Minkowski geometry were developed solely using this system of numbers ( see, \cite{Ca-Ca-Ca-Za03},\cite{Ca-Ca-Ca-Za04},\cite{Fj86},\cite{So95},\cite{Ya79}). During the past several years research in this area has been to develop hyperbolic numbers as an affordable replacement for the real number system. Many papers has appeared studying hyperbolic numbers from various points of view (see, \cite{Ca-Bo-Ca08}, \cite{Mo-Ro-98},  
\cite{Ca-Za11}, \cite{Gu-Sp97} ,\cite{He-So97}) and references therein. However in recent paper \cite{Ga-Ko16}  studied  this system of numbers as the only(natural) generalization of real numbers, into Archimedean f-algebra of dimension two. They have also generalized the fundamental properties of real numbers to this number system. The set of bicomplex numbers and the hyperbolic number system seems to have originated independently. Recently, a lots of work is being done on  bicomplex numbers and bicomplex functional analysis. However later it was found that hyperbolic numbers is a subset of the set of bicomplex numbers and it plays the same role for bicomplex numbers as real numbers plays for the set of complex numbers.\\
Bicomplex numbers are being studied for quite a long time now. A lot of work have already been published in recent years. The book of G. B. Price \cite{Pr91} contains a comprehensive study of bicomplex numbers. Now the theory of functional analysis with bicomplex scalars is essentially a new subject of study. The study of bicomplex functional analysis was initiated in \cite{GE-MA-RO1} and \cite{GE-MA-RO3}. The monograph \cite{Al-Lu-Sh-Sr13} is a systematic introductory text on functional analysis with bicomplex scalars. For more details on bicomplex numbers and bicomplex functional analysis  ( see,\cite{Co-Sa-St11},\cite{De-Va-Va12},\cite{Ku-Sa16},\cite{Ku-Si15},\cite{Lu-Sh-Va12} \cite{Pr91})and (\cite{GE-MA-RO1}, \cite{GE-MA-RO2}, \cite{GE-MA-RO3}) and references therein.\\

The concept of convexity plays a vital role in the study of analysis. One such concept is that of the study of convex function which is specially important in the study of optimization problems. In the study of fourier analysis W.H. Young found certain convex functions $\Phi : \mathbb{R} \to \bar{\mathbb{R}}^+$ which satisfies $\Phi(0)=0$, $\Phi(-x)=\Phi(x)$ and $\lim_{x \to \infty}\Phi(x)=\infty$ The importance of Young functions were only recognized from the work of W.H. Young in $1912$. But their role in the abstract analysis emerges only with the fundamental researches of Z.W. Birnbaum and W.Orlicz in the year 1931. Here we shall define the same types of convex function from the set of hyperbolic number $\mathbb{D}$ to the set of hyperbolic number $\mathbb{D}$. Further we shall define hyperbolic valued Young function. For details on Orlicz spaces, we refer to \cite{Rao-Ren91} and  \cite{Rao-Ren02}. Our results on Orlicz spaces with bicomplex scalars are essentially based on the proofs of corresponding results on Orlicz spaces with complex scalars as in \cite{Rao-Ren91}. 
\end{section}

\section{\bf A review of bicomplex and hyperbolic numbers}
\medskip
First we shall summarize briefly some basic properties about bicomplex numbers, bicomplex holomorphic functions and bicomplex linear operators. As in case of complex numbers, we take only one imaginary unit, but for bicomplex numbers, we consider two imaginary units $i$ and $j$ satisfying $ij=ji$ with $i^2=j^2=-1.$ Let $\mathbb{C}(i)$ be the set of complex numbers with imaginary units $i$ and let $\mathbb{C}(j)$ be the set of complex numbers with imaginary units $j$.  Then we define set of bicomplex numbers denoted by $\mathbb{BC}$ as 
\begin{equation}
\label{eq:F1}
\mathbb{BC}=\left\{Z=x_o+x_1i+x_2j+x_3ij: x_o,x_1,x_2,x_3 \in \mathbb{R}\right\}= \left\{Z=z_1+jz_2: z_1,z_2 \in \mathbb{C}(i)\right\}
\end{equation}
The set $\mathbb{BC}$ becomes a ring under the operation of addition and multiplication defined by
$$Z+W=(z_1+jz_2)+(w_1+jw_2)=(z_1+w_1)+j(z_2+w_2)$$
$$Z\cdot W==(z_1+jz_2)(w_1+jw_2)=(z_1w_1-z_2w_2)+j(z_2w_1+z_1w_2)$$
and also $\mathbb{BC}$ becomes a module over itself. The set of complex numbers $\mathbb{C}(i)$ is a subring of $\mathbb{BC}.$ Due to the fact that $\mathbb{BC}$ has two imaginary unit whose square is $-1$ and one hyperbolic unit whose square is $1$ . There are three conjugations for $\mathbb{BC}$.
\begin{itemize}
\item[(i)]The bar-conjugation\\
$\bar{Z}=\bar{z_1}+j \bar{z_2}$\\
\item[(ii)]The $\dagger$-conjugation\\
$Z^{\dagger}=z_1-jz_2$\\

\item[(iii)] The $\ast$-conjugation\\
$Z^\ast=\bar{Z}^{\dagger}=\bar{z_1}-j\bar{z_2},$\\
\end{itemize}
where $\bar{z_1}$,$\bar{z_2}$ are the usual conjugation of complex numbers $z_1$,$z_2$ in $\mathbb{C}(i).$\\
Each of the conjugation is an additive, multiplicative and involutive operation on $\mathbb{BC}$ and is a ring automorphism of $\mathbb{BC}.$ Three possible moduli arise in accordance to the multiplication of a bicomplex number by its three different conjugations.
\begin{itemize}
\item[(i)] $ Z \cdot Z^{\dagger}={|Z|_i}^2={z_1}^2+{z_2}^2=(|\eta_1|^2-|\eta_2|^2)+2Re(\eta_1 \eta_2^{\ast})i.$\\
\item[(ii)] $ Z \cdot \bar{Z}={|Z|_j}^2={\eta_1}^2+{\eta_2}^2=(|z_1|^2-|z_2|^2)+2Re(z_1 \bar{z_2})j.$\\
\item[(iii)] $ Z \cdot Z^{\ast}={|Z|_k}^2=(|z_1|^2+|z_2|^2)-2Im(z_1 \bar{z_2})k.$\\
\end{itemize}

It is to be noted that unlike in the complex case, each of these modules are not $\mathbb{R}^+$-valued, but $\mathbb{C}(i)$, $\mathbb{C}(j)$ and $\mathbb{D}$-valued functions. Now the $\dagger$-conjugation is particularly important as it paved  us way to define the invertiblity of a bicomplex number. Any bicomplex number $Z$ is said to be invertible if $Z \cdot Z^{\dagger} \ne 0$ and its inverse is given by
$$Z^{-1}=\frac{Z^{\dagger}}{Z \cdot Z^{\dagger}}=\frac{Z^{\dagger}}{|Z|_i^2}.$$
Also, if $Z \ne 0$, but $Z \cdot Z^{\dagger}=|Z|_i=0$ , then $Z$ is said to be a zero-divisor. We denote the set of all zero-divisors by
$$\mathcal{NC}=\left\{Z=z_1+jz_2: Z \ne 0,Z \cdot Z^{\dagger}=z_1^2+z_2^2=0\right\}$$
and is called the Null cone of $\mathbb{BC}$. Now there are two very special zero divisors. We call them idempotent elements defined as
$$e_1=\frac{1}{2}(1+ij)~~~\mbox{and}~~~e_2=\frac{1}{2}(1-ij)$$ where $z_1=\frac{1}{2}$ and $z_2=i\frac{1}{2},$ considering $e_1$ and $e_2$ as bicomplex numbers although these are hyperbolic numbers which is a subring of $\mathbb{BC}.$
It is easy to see that $e_1$ and $e_2$ are zero divisors and are mutually complementary idempotent elements.
The sets $\mathbb{BC}_{e_1}=e_1\mathbb{BC}$ and $\mathbb{BC}_{e_2}=e_2\mathbb{BC}$ are (principal) ideals in the ring $\mathbb{BC}$ and have the property that
$$\mathbb{BC}_{e_1}\cap \mathbb{BC}_{e_2}=\{0\},$$ and
\begin{equation}
\label{eq:F3}
\mathbb{BC}=\mathbb{BC}_{e_1}+\mathbb{BC}_{e_2}.
\end{equation}
Equation \eqref{eq:F3} is sometimes called the idempotent decomposition of the ring of bicomplex numbers $\mathbb{BC}.$
The Euclidean norm $|\cdot|$ of a bicomplex number $Z=z_1e_1+z_2e_2$ is defined  as $|Z|=\sqrt{{{x_1}^2+{x_2}^2+{x_3}^2+{x_4}^2}}=\sqrt{{|z_1|}^2+{|z_2|}^2}$ and one can easily check that for any $Z$ and $W$ in $\mathbb{BC}$, we have 
$$|Z \cdot W| \leq \sqrt{2}|Z||W|.$$
The $\mathbb{D}$-valued norm of the bicomplex number $Z=z_1e_1+z_2e_2$ denoted  by $|Z|_k$ is defined as $|Z|_k=|z_1|e_1+|z_2|e_2$ where $|z_1|$ and $|z_2|$ are the  usual modulus of complex numbers $z_1$ and $z_2$. Further $|Z\cdot W|_k=|Z\cdot|_k|W\cdot|_k,$ and Euclidean norm and hyperbolic norm  of a bicomplex number  is related by $$\left||Z|_k\right|=\left|Z\right|.$$ For more details on Euclidean norm and hyperbolic norm one can refer to \cite{Al-Lu-Sh-Sr13}. It is interesting to mention here that any $\mathbb{BC}$-module X has some similar structural properties of the scalars ring $\mathbb{BC}.$ For instance we can also decompose the $\mathbb{BC}$-module X as follows.
\begin{equation}
\label{eq:R1}
X=X_1e_1+X_2e_2
\end{equation}
where $X_1=e_1X$ and $X_2=e_2X$ are $\mathbb{C}(i)$-linear  space as well as a $\mathbb{BC}$ module. Equation $\eqref{eq:R1}$ is called the idempotent representation of the $\mathbb{BC}$-module X. Thus each $x \in X$  can uniquely be expressed as $x=x_1e_1+x_2e_2$ and it allows us with component wise addition and scalar multiplication of elements in $\mathbb{BC}$-module X and scalars  in $\mathbb{BC}.$ A real valued norm on $\mathbb{BC}$-module X  is defined as 
\begin{equation}
\label{eq:RNorm}
\left\|x\right\|=\sqrt{\frac{\|x_1\|_1^2+\|x_2\|_1^2}{2}}
\end{equation}
for any $x \in X,$  where $\|.\|_1$ and $\|.\|_2$ are real valued norms on $X_1$ and $X_2.$ However for any scalar $\alpha \in \mathbb{BC}$ and $x \in X$, the norm satisfies the inequality $\|\alpha x\| \leq \sqrt{2} |\alpha| \|x\|.$  Although the ring $\mathbb{BC}$ contains three two-dimensional subrings which are also forming real subalgebras $\mathbb{C}(i)$, $\mathbb{C}(j)$ and $\mathbb{D}$. However first two is isomorphic to $\mathbb{C}$, and so it is not wise to study them separately in the context of the set of bicomplex numbers.  But the study of the set of hyperbolic number $\mathbb{D}$ is interesting and important in itself as we can see in the recent paper by Hichem Gargoubi and Sayed Kossentini \cite{Ga-Ko16}. 
\\
\\
Now we shall also go through a brief review of  hyperbolic numbers system.
The set of hyperbolic numbers denoted by $\mathbb{D}$ is a ring of all numbers of the form $Z=a+bk,$ a,b $\in  \mathbb{R}$, with $k$ satisfying $k^2=1.$
$$i.e~~~~~~~~~~~~~~ \mathbb{D}=\left\{a+bk: a,b \in \mathbb{R}, k^2=1, k \notin \mathbb{R}\right\}.$$
The set $\left\{x+yij:x,y \in \mathbb{R},i^2=j^2=-1\right\}$ is a subset of the set of bicomplex numbers which is isomorphic to $\mathbb{D}$ as a real algebra. It is very well known that we can also decompose the set of hyperbolic as
\begin{equation}
\label{eq:1d}
\mathbb{D}=\mathbb{D}e_{1}+\mathbb{D}e_{2}
\end{equation}
 where $e_{1}=\frac{1}{2}(1+ij)$ and  $e_{2}=\frac{1}{2}(1-ij)$ are two hyperbolic numbers belonging to the null cone. We call equation \eqref{eq:1d}, the idempotent decomposition of $\mathbb{D}$.  Thus the idempotent representation of any hyperbolic number $x=x_1+x_2k$ is given by 
$$x=\alpha e_{1}  + \beta e_{2},~~~ \alpha, \beta \in \mathbb{R}$$
with $\alpha=x_2+x_1$, $\beta=x_2-x_1$. Also the set of positive hyperbolic number denoted by $\mathbb{D}^+$ are the set of all those hyperbolic numbers whose idempotent components are non negative.
$$i.e~~~~\mathbb{D}^+=\left\{\alpha e_{1} +\beta e_{2}:\alpha, \beta \geq 0\right\}.$$ We shall define a partial order relation on $\mathbb{D}$ as follows.\\ 
Given $x,y \in \mathbb{D},$ we write $x \leq' y$ if $y-x \in \mathbb{D}^+$. It is easy to see that this relation is reflexive, symmetric and antisymmetric and so it defines a partial order relation on $\mathbb{D}.$  Also for $x ,y \in \mathbb{D}$, if $x\leq' y$, then we say that $y$ is $\mathbb{D}$-larger than $x$ and $x$ is $\mathbb{D}$-smaller than $y$. The notion of upper and lower bounds also exists in the context of hyperbolic plane. Given a subset $S$ of $\mathbb{D}$ we can define $\mathbb{D}$-upper bounds and $\mathbb{D}$-lower bounds of this set $S$. Using this bounds, this set can be made $\mathbb{D}$-bounded from above and $\mathbb{D}$-bounded from below if it exists. Now if the set is $\mathbb{D}$-bounded from above as well as from below then we say that the set is $\mathbb{D}$-bounded.\\
We can also define hyperbolic-valued($\mathbb{D}$-valued) norm  on $\mathbb{BC}$-module X as follows:
\begin{equation}
\label{eq:NHy}
\|x\|_{\mathbb{D}}=\|x_1e_1+x_2e_2\|_{\mathbb{D}}= \|x_1\|_1e_1+\|x_2\|_2e_2.
\end{equation}
such that for any $x,y \in X$ and scalar $\alpha \in \mathbb{BC}$ we have
\begin{equation}
 \|\alpha x\|_{\mathbb{D}} = |\alpha|_k \cdot \|x\|_{\mathbb{D}}~~~\mbox{and}~~~\|x+y\|_{\mathbb{D}} \leq ' \|x\|_{\mathbb{D}}+\|y\|_{\mathbb{D}}.
\end{equation}
where $\leq '$ is a partial order relation on $\mathbb{D}$.\\
We further define the notion of the supremum of a given subset of $\mathbb{D}$. Supremum of $S \subset \mathbb{D}$ denoted by $\sup_{\mathbb{D}}S$ is defined usually as the least of all $\mathbb{D}$-upper bounds of the given set. Similarly $\inf_{\mathbb{D}}S$ is the greatest of all $\mathbb{D}$-lower bounds of the set. However due to the idempotent decomposition of $\mathbb{D}$, we can find a convenient expression as follows:\\
For a subset $S$ of $\mathbb{D}$ which is $\mathbb{D}$- bounded from above, we consider the set $C_1=\left\{\alpha:\alpha e_{1}+ \beta e_{2} \in S\right\}$ and $C_2=\left\{\beta:\alpha e_{1}+ \beta e_{2} \in S\right\}.$ Then the supremum of the set S denoted by $\sup_{\mathbb{D}}S$ is defined as
$\sup_{\mathbb{D}}S=\sup C_1 \cdot e_{1} + \sup C_2 \cdot e_{2}$ where $\sup C_1$ and $\sup C_2$ is the supremum taken over the subset $C_1$ and $C_2$ of real numbers. Finally the hyperbolic modulus of any hyperbolic number $Z=\alpha e_{1} + \beta e_{2}$ denoted by $\left|Z\right|_k$ is given by the formula $\left|Z\right|_k=Z \cdot Z^*= \left|\alpha\right|^2 e_{1} + \left|\beta\right|^2 e_{2}$ where "*" denotes the ${*}$-conjugation. For further basic results on hyperbolic numbers (see \cite{Al-Lu-Sh-Sr13} \cite{Ga-Ko16} and \cite{LUN-REgaldo}).\\

The paper is organized as follows. In Section 2, we discuss some definitions including the convexity of a function from $\mathbb{D}$ to $\mathbb{D}$. In this section we have also obtained the integral representation of a convex function from $\mathbb{D}$ to $\mathbb{D}.$ Further we have the hyperbolic version of Young's inequality for a pair of Young functions(N-functions). Section 3 consists of a short introduction to the $\mathbb{D}$-metric spaces and we have discussed conditions for the completeness of a $\mathbb{D}$-metric spaces. In section 4, we have discussed the main topic of this paper which is Orlicz Space with bicomplex scalars. Here we have studied some properties of this space  equipped with the hyperbolic valued Luxemburg norm and also proved some results on it. In section 5, we have discussed the spectrum of bicomplex linear operators. Since from the very well known theory of spectrum of linear operators, we know that the spectrum of complex  linear operator is bounded, however the spectrum of bicomplex linear operator is unbounded \cite{Co-Sa-St13}. Here we have proved that that the point spectrum of unilateral shift operator on $l^2(\mathbb{BC})$ is in the null cone. Finally in the last section, we define Musielak Orlicz spaces with bicomplex scalars.
\\ 
\\
\\

{\bf Preliminaries} 
We shall begin with some definition which will be used throughout this paper.
\begin{definition}
Let $\Omega_{\mathbb{D}}$ denotes a non-empty set and $\Sigma_{\mathbb{D}}$ be the sigma algebra of subsets of $\Omega_{\mathbb{D}}$. Then by a hyperbolic measure on $\Sigma_{\mathbb{D}}$, we  mean a set function $\mu_{\mathbb{D}}: \Sigma_{\mathbb{D}} \to {\mathbb{\overline{D}}}_+$
\\
\begin{itemize}
\item[(i)] $\mu_{\mathbb{D}}$ is non-negative extended hyperbolic valued set function.\\
\item[(ii)] $\mu_{\mathbb{D}}(\phi)=0$ and\\
\item[(iii)] $\mu_{\mathbb{D}}$ is countable additive, i.e., for any disjoint subcollection $\{E_{n}:n \in \mathbb{N}\} \subset \Sigma_{\mathbb{D}},$
\end{itemize}
$$\mu_{\mathbb{D}}\left(\bigcup_{n \in \mathbb{N}}E_n\right)= \sum_{n \in \mathbb{N}}\mu_{\mathbb{D}}\left(E_n\right).$$ 
Every $\mathbb{D}$-valued measure can be written as
$\mu_{\mathbb{D}}(E)=\psi_1(E)+k\psi_2(E) =\mu_{\mathbb{D}}^1(E) e_1+ \mu_{\mathbb{D}}^2(E) e_2$ with $\mu_{\mathbb{D}}^1(E) = \psi_1(E)+\psi_2(E)$ and  $\mu_{\mathbb{D}}^2(E) = \psi_2(E)-\psi_1(E)$ for every $E \in \Sigma_{\mathbb{D}}.$ Now by property (i) of $\mathbb{D}$-valued measure implies that
$$\mu_{\mathbb{D}}^1(E) \geq 0 ~~\mbox{and}~~ \mu_{\mathbb{D}}^2(E) \geq 0~~\mbox{for every}~ E \in \Sigma_{\mathbb{D}}.$$

Further property (iii) of $\mu_{\mathbb{D}}$ implies that
$$\mu_{\mathbb{D}}^1\left(\bigcup_{n \in \mathbb{N}}E_n\right) e_1+ \mu_{\mathbb{D}}^2\left(\bigcup_{n \in \mathbb{N}}E_n\right) e_2=\sum_{n \in \mathbb{N}} \mu_{\mathbb{D}}^1(E_n) e_1+ \sum_{n \in \mathbb{N}} \mu_{\mathbb{D}}^2(E_n) e_2$$ 
Thus $\mu_{\mathbb{D}}$ is called a $\mathbb{D}$-valued measure on the sigma algebra $\Sigma_{\mathbb{D}}$ and the triplet $\left(\Omega_{\mathbb{D}},\Sigma_{\mathbb{D}},\mu_{\mathbb{D}}\right)$ is called a $\mathbb{D}$-measure spaces.
\end{definition}
\begin{remark}
By non negative extended hyperbolic number we mean the set
$$\overline{\mathbb{D}}^{+}=\left\{ae_1+be_2:a,b>0\right\}\cup\left\{\infty\right\}\cup\{-\infty\}\cup\{\infty e_1+be_2\}\cup \{ae_1-\infty e_2\}$$
\end{remark} 
\begin{example} {(\cite{Ap-Lu-Sh15})}
Let $\left(\Omega_{\mathbb{D}},\Sigma_{\mathbb{D}}\right)$ be a measurable space, a $\mathbb{D}$-valued function
$$\mu_{\mathbb{D}}: \Sigma_{\mathbb{D}} \longmapsto \mathbb{D}$$
with the property\\
\begin{itemize}
\item[(i)] $\mu_{\mathbb{D}}(E) \geq ' 0$ for every $E \in \Sigma_\mathbb{D}.$\\
\item[(ii)] $\mu_{\mathbb{D}}(\Omega_{\mathbb{D}})=\delta$, where $\delta$ takes only three possible values $1,e_1,e_2.$\\
\item[(iii)] $\mu_{\mathbb{D}}$ satisfies the countable additive property.\\
\end{itemize}
is a $\mathbb{D}$-valued probabilistic measure or $\mathbb{D}$-valued probability, on the sigma algebra of all events in $\Sigma_{\mathbb{D}}$.
\end{example}

\begin{definition}
A $\mathbb{D}$-valued measure $\mu_{\mathbb{D}}$ on the sigma algebra $\Sigma_{\mathbb{D}}$ of subsets of $\Omega_{\mathbb{D}}$ is said to be a $\mathbb{D}$-finite measure if $\mu_{\mathbb{D}} \in \mathbb{D}^+$ but $\mu_{\mathbb{D}} \notin \mathbb{\overline{D}}_+.$
In this case we say that $\(\Omega_{\mathbb{D}},\Sigma_{\mathbb{D}},\mu_{\mathbb{D}}\)$ is a $\mathbb{D}$-finite measure space.
\end{definition}
\begin{definition}
A sigma algebra $\Sigma_{\mathbb{D}}$ of subsets of $\Omega_{\mathbb{D}}$ is said to be complete with respect to the $\mathbb{D}$-valued measure $\mu_{\mathbb{D}}$, if for every $A_o \subset A$ with $\mu_{\mathbb{D}}(A)=0$ implies that $A_0 \in \Sigma_{\mathbb{D}}$ and we say that 
$\left(\Omega_{\mathbb{D}},\Sigma_{\mathbb{D}},\mu_{\mathbb{D}}\right)$ is a $\mathbb{D}$-complete measure space.
\end{definition}

\begin{definition}
Let $\left(\Omega_{\mathbb{D}}, \Sigma_{\mathbb{D}}, \mu_{\mathbb{D}}\right)$ be a $\mathbb{D}$-measure space. A $\mathbb{BC}$-valued function $f$ defined on  a subset $S$ of $\Omega_{\mathbb{D}}$ is said to be a $\Sigma_{\mathbb{D}}$- measurable function, if the subset of $S$ consisting of all those elements $t \in S$ for which $\alpha-|f(t)|_k \in \mathbb{D}^+$ is an element of $\Sigma_{\mathbb{D}}$ for every $\alpha \in \mathbb{D}.$
$$i.e.,~~~~~~~~~~\left\{t \in S: |f(t)|_k \leq' \alpha \right\} \in  \Sigma_{\mathbb{D}} .$$
Further, since that as $f$ is a $\mathbb{BC}$-valued $\Sigma_{\mathbb{D}}$-measurable function and so we can decompose $f$ as 
$$|f|_k=|f_1| e_1+|f_2| e_2$$ where $f_1=e_1 f$ and $f_2=e_2 f$ are complex valued $\Sigma_{\mathbb{D}}$-measurable functions on the measure spaces $\left(\Omega_{\mathbb{D}}, \Sigma_{\mathbb{D}}, \mu_{\mathbb{D}}^1\right)$ and $\left(\Omega_{\mathbb{D}}, \Sigma_{\mathbb{D}}, \mu_{\mathbb{D}}^2\right)$ respectively where  $\mu_{\mathbb{D}}^1=e_1\mu_{\mathbb{D}}$ and $\mu_{\mathbb{D}}^2=e_2\mu_{\mathbb{D}}$ are real valued measures on $\Omega_{\mathbb{D}}$.
\end{definition}
\begin{definition}(\cite{Al-Lu-Sh-Sr13}){($\mathbb{D}$-Derivable functions)}
Let $f:A \subset \mathbb{BC} \to \mathbb{BC}$ be a function. Then the derivative of $f$ at a point $x_o \in A$ is defined to be the limit $f'(x_o)$, if it exists, and is given by 
$$f'(x_o)= \lim_{h \to 0}\frac{f(x_o+h)-f(x_o)}{h}$$
for every $\mathbb{D}$-small bicomplex number $h\in \mathbb{BC} \setminus \mathcal{NC}.$ 
Now if $f'(x)$ exists for every $x \in A$, then we say that $f$ is derivable on A. We can obtain the second order derivative and so on if it exists.
\end{definition}
\begin{definition}
A $\mathbb{BC}$-valued function for which derivative of all order exists and are continuous is call  $C^{\infty}$ function. We define  $C^{\infty}_C(\mathbb{BC})$ to be the collection of all $C^{\infty}$ functions with compact support.
$$i.e.,~~~C^{\infty}_C(\mathbb{BC})=\{f: f~~\mbox{is}~~C^{\infty}~~ \mbox{function with}~\overline{Supp(f)}~~ \mbox{compact}\}$$

\end{definition}

\begin{definition}(\cite{Mo-Ro-98}) Let $f: \mathbb{D} \to \mathbb{D}$ given by  $f(h)=u(x,t)+k v(x,t)$ be a hyperbolic continuous function where $h=x+kt$ is a hyperbolic number, and $\gamma$ is a Jordan arc with $\gamma :[a,b] \to \mathbb{D}$ and $\gamma^{'}(t)$ is a continuous except for a finite number of points. Then the hyperbolic integral of $f$ on $\gamma$ is define as the line integral 
\begin{equation}
\label{eq:LI1}
\int_{\gamma}f(h)dh= \int_{\gamma}[u(h)dx+v(h)dt]+k \int_{\gamma}[u(h)dt+v(h)dx]
\end{equation}
\end{definition}

We can decompose the hyperbolic integral as follows.
\begin{remark}
Since $f: \mathbb{D} \to \mathbb{D}$ implies that  $f(h)=u(x,t)+k v(x,t)=f_1(x,t)e_1+f_2(x,t)e_2$, where $f_1(x,t)=u(x,t)+v(x,t)$ and $f_1(x,t)=v(x,t)-u(x,t)$. Since the integral defined in equation \eqref{eq:LI1} is hyperbolic valued and so from equation \eqref{eq:LI1} we have
\begin{eqnarray*}
\int_{\gamma}f(h)dh & = & \int_{\gamma}[u(h)dx+v(h)dt]+k \int_{\gamma}[u(h)dt+v(h)dx]\\
                    & = & \left\{ \int_{\gamma}[u(h)dx+v(h)dt] + \int_{\gamma}[u(h)dt+v(h)dx \right\}e_1 \\
										& + & \left\{\int_{\gamma}[u(h)dt+v(h)dx -\int_{\gamma}[u(h)dx+v(h)dt]  \right\}e_2\\						
\end{eqnarray*}
Now rearranging the terms and simplifying, we get

\begin{eqnarray*}
\int_{\gamma}f(h)dh & = & \left\{ \left(\int_{\gamma}u(h)dx+\int_{\gamma}v(h)dx\right)+k\left(\int_{\gamma}u(h)dt+\int_{\gamma}v(h)dt\right)\right\}e_1\\
&-&	\left\{\left(\int_{\gamma}v(h)dx-\int_{\gamma}u(h)dx\right)+k\left(\int_{\gamma}v(h)dt+\int_{\gamma}u(h)dt\right)\right\}e_2	\\
&=& \left\{ \left(\int_{\gamma}u(h)dx+k\int_{\gamma}u(h)dt+\int_{\gamma}v(h)dx+k\int_{\gamma}v(h)dt\right)\right\}e_1\\
&-&	\left\{\left(\int_{\gamma}v(h)dx+k\int_{\gamma}v(h)dt-\int_{\gamma}u(h)dx-k\int_{\gamma}u(h)dt\right)\right\}e_2\\
&=& \left\{ \left(\int_{\gamma}u(h)d(x+kt)+\int_{\gamma}v(h)d(x+kt)\right)\right\}e_1\\
&-&	\left\{\left(\int_{\gamma}v(h)d(x+kt)-\int_{\gamma}u(h)d(x+kt)\right)\right\}e_2\\
&=& \left\{ \left(\int_{\gamma}(u(h)+v(h))d(x+kt)\right)\right\}e_1 -	\left\{\left(\int_{\gamma}(v(h)-u(h))d(x+kt)\right)\right\}e_2\\
&=& \left\{ \left(\int_{\gamma}(u(h)+v(h))d(h)\right)\right\}e_1 -	\left\{\left(\int_{\gamma}(v(h)-u(h))d(h)\right)\right\}e_2\\
&=& \left\{\int_{\gamma}f_1(h)d(h)\right\}e_1 -\left\{\int_{\gamma}f_2(h)d(h)\right\}e_2\\
\end{eqnarray*}
Thus, we get the idempotent decomposition of the hyperbolic integral as
\begin{equation}
\int_{\gamma}f(h)dh = \left\{\int_{\gamma}f_1(h)dh\right\}e_1 -\left\{\int_{\gamma}f_2(h)dh\right\}e_2.
\end{equation}
where $f(h)=u(h)+k v(h)=f_1(h)e_1+f_2(h)e_2$ and $f_1(h)=u(h)+v(h)$ and $f_2(h)=v(h)-u(h)$
\end{remark}

Now we shall define $\mathbb{D}$-valued convex function.
\begin{definition} {$\(\mathbb{D}\mbox{-Valued Convex Function}\)$} A function $\varphi_{\mathbb{D}}:\mathbb{D} \to \mathbb{\overline{D}}_+$ is said to be a $\mathbb{D}$-valued convex function if for every $x,y \in \mathbb{D}$ with $0 \leq' \alpha \leq' 1$ , we have that
$$\varphi_{\mathbb{D}}(\alpha x + (1-\alpha) y) \leq' \alpha \varphi_{\mathbb{D}}(x)+ (1-\alpha)\varphi_{\mathbb{D}}(y).$$
\end{definition}
\medskip
\begin{example}
Let $\varphi_{\mathbb{D}} : \mathbb{D} \to \mathbb{D}^+$ be defined by $\varphi_{\mathbb{D}}(x)=|x|_k$ for every $x \in \mathbb{D}$. Then clearly $\varphi_{\mathbb{D}}$ is a $\mathbb{D}$-convex function.  
\end{example}

Note that if $\varphi_{\mathbb{D}}:\mathbb{D} \to \mathbb{D}^+$ is a function  such that it can be written as 
$$\varphi_{\mathbb{D}}(x)=\psi_1(x) + \psi_2(x) k = \varphi_{\mathbb{D}_1}(x)e_{1} + \varphi_{\mathbb{D}_2}(x) e_{2},$$
 where $\psi_1, \psi_2$ are functions from $\mathbb{D} \to \mathbb{R}$ satisfying $(\psi_1)^2 - (\psi_2)^2 \geq 0$ and $(\psi_1) \geq 0 $ and $\varphi_{\mathbb{D}_1},\varphi_{\mathbb{D}_1}$ are  the idempotent components  given by 
$$\varphi_{\mathbb{D}_1}=\psi_2+\psi_1 \geq 0~~ \mbox{and}~~ \varphi_{\mathbb{D}_2}=\psi_2-\psi_1 \geq 0 .$$
 We further define the decomposition of $\varphi_{\mathbb{D}}$ as
\begin{equation} \label{eq:1}
\varphi_{\mathbb{D}}(x)=\varphi_{\mathbb{D}_1}(e_1x)e_{1} + \varphi_{\mathbb{D}_2}(e_2x) e_{2}.
\end{equation}

Next we shall show that if $\varphi_{\mathbb{D}}$ is a $\mathbb{D}$-valued convex function, then the corresponding idempotent components $\varphi_{\mathbb{D}_1}$ and $\varphi_{\mathbb{D}_2}$ are $\mathbb{R}$-convex functions.

\begin{theorem}
Let $\varphi_{\mathbb{D}}: \mathbb{D} \to \mathbb{D}^+$ be a $\mathbb{D}$-convex function, and let
$$\varphi_{\mathbb{D}}(x)=\varphi_{\mathbb{D}_1}(e_{1}x) e_{1} + \varphi_{\mathbb{D}_2}(e_{2}x) e_{2}$$ be the idempotent representation of $\varphi_{\mathbb{D}}$ where $\varphi_{\mathbb{D}_1}:e_1\mathbb{D} \to \mathbb{R}$ and $\varphi_{\mathbb{D}_2}:e_2\mathbb{D} \to \mathbb{R}$. Then $\varphi_{\mathbb{D}_1}$ and $\varphi_{\mathbb{D}_2}$ are $\mathbb{R}$-convex functions.  
\end{theorem}
\begin{proof}
Suppose that $\varphi_{\mathbb{D}}$ is a $\mathbb{D}$-convex function. We shall first show that $\varphi_{\mathbb{D}_1}$ is $\mathbb{R}$-convex function. For this let $x',y' \in e_1\mathbb{D}$, so that we can find $x,y \in \mathbb{D}$ such that $x'=e_1x$ and $y'=e_1y.$ Suppose that $0 \leq \lambda_1 \leq 1$ and $0 \leq \lambda_2 \leq 1$ be such that
$0 \leq'  \lambda \leq' 1,$ where $\lambda=\lambda_1 e_1+\lambda_2 e_2.$ Now since $\varphi_{\mathbb{D}}$ is a $\mathbb{D}$-convex function, 
$$\varphi_{\mathbb{D}}(\lambda x + (1-\lambda) y) \leq' \lambda \varphi_{\mathbb{D}}(x)+ (1-\lambda)\varphi_{\mathbb{D}}(y).$$
This further implies that
\begin{eqnarray*}
\varphi_{\mathbb{D}_1}\left(e_1(\lambda x+(1-\lambda)y)\right)e_1 + \varphi_{\mathbb{D}_2}\left(e_2(\lambda x+(1-\lambda)y)\right)e_2 &\\
 \leq' \lambda \left(\varphi_{\mathbb{D}_1}(e_1 x)e_1+ \varphi_{\mathbb{D}_2}(e_2 x)e_2\right) & \\
 + (1-\lambda) \left(\varphi_{\mathbb{D}_1}(e_1 y)e_1+ \varphi_{\mathbb{D}_2}(e_2 y)e_2\right). &
\end{eqnarray*}
Now multiplying the above equation by $e_1$, we get
\begin{equation}
\label{eq:N1}
\varphi_{\mathbb{D}_1}\left(e_1(\lambda x + (1-\lambda) y)\right)e_1 \leq' \left(\lambda \varphi_{\mathbb{D}_1}(e_1 x)+(1-\lambda) \varphi_{\mathbb{D}_1}(e_1 y)\right)e_1
\end{equation}
From LHS of \eqref{eq:N1}, we get

\begin{eqnarray*}
\varphi_{\mathbb{D}_1}\left(e_1(\lambda x + (1-\lambda) y)\right)e_1 & = & \varphi_{\mathbb{D}_1}\left(e_1((\lambda_1e_1+\lambda_2 e_2) x + ((1-\lambda_1)e_1+ (1-\lambda_2)e_2 )y)\right)e_1 \\
& = & \varphi_{\mathbb{D}_1}\left(e_1 \lambda_1 x + e_1(1-\lambda_1)y\right)e_1\\
& = & \varphi_{\mathbb{D}_1}\left(\lambda_1 x' + (1-\lambda_1)y'\right)e_1.\\
\end{eqnarray*}
\begin{equation}
\label{eq:N2}
\mbox{So,}~~~\varphi_{\mathbb{D}_1}\left(e_1(\lambda x + (1-\lambda) y)\right)e_1  =  \varphi_{\mathbb{D}_1}\left(\lambda_1 x' + (1-\lambda_1)y'\right)e_1
\end{equation}
Similarly solving RHS of \eqref{eq:N1}, we get
\begin{equation}
\label{eq:N3}
\left(\lambda \varphi_{\mathbb{D}_1}(e_1 x)+(1-\lambda) \varphi_{\mathbb{D}_1}(e_1 y)\right)e_1  =  \left(\lambda_1 \varphi_{\mathbb{D}_1}(x')+(1-\lambda_1)\varphi_{\mathbb{D}_1}(y')\right)e_1.
\end{equation}

Using equation \eqref{eq:N2} and \eqref{eq:N2} in $\eqref{eq:N1}$, we finally get

\begin{eqnarray}
\varphi_{\mathbb{D}_1}\left(\lambda_1 x' + (1-\lambda_1)y'\right)e_1 & \leq' & \left(\lambda_1 \varphi_{\mathbb{D}_1}(x')+(1-\lambda_1)\varphi_{\mathbb{D}_1}(y')\right)e_1
\end{eqnarray}
Thus,
$$\varphi_{\mathbb{D}_1}\left(\lambda_1 x' + (1-\lambda_1)y'\right) \leq \lambda_1 \varphi_{\mathbb{D}_1}(x')+(1-\lambda_1)\varphi_{\mathbb{D}_1}(y')$$
This proves that $\varphi_{\mathbb{D}_1}$ is a $\mathbb{R}$-valued convex function. Similarly we can prove for $\varphi_{\mathbb{D}_2}.$ 
This completes the proof.
\end{proof}

\begin{example}
Let $\varphi_{\mathbb{D}_1}:e_1\mathbb{D} \to \mathbb{R}$ and $\varphi_{\mathbb{D}_2}:e_2\mathbb{D} \to \mathbb{R}$   be defined by 
$$\varphi_{\mathbb{D}_1}(x)=a(x+e_2\mathcal{P}_1(x))~~\mbox{and}~~\varphi_{\mathbb{D}_2}(x)=a(x+e_1\mathcal{P}_2(x)),$$
where $\mathcal{P}_1$ and $\mathcal{P}_2$ are projections onto the first and second coordinate from $\mathbb{D}$ to $e_1\mathbb{R}+e_2 \mathbb{R}$. Then clearly $\varphi_{\mathbb{D}_1}$ and $\varphi_{\mathbb{D}_2}$ are $\mathbb{R}$-convex functions. Now if we take $\varphi_{\mathbb{D}}(x)=\varphi_{\mathbb{D}_1}(e_1x)e_1+\varphi_{\mathbb{D}_2}(e_2x)e_2$, then it is easy to show that $\varphi_{\mathbb{D}}$ is a $\mathbb{D}$-valued convex function.
\end{example}

We shall now prove this for any two $\mathbb{R}$-convex functions.
\begin{theorem}
Let $\varphi_1 : e_1 \mathbb{D} \to \mathbb{R}^+$ and $\varphi_2 : e_2 \mathbb{D} \to \mathbb{R}^+$ be any two $\mathbb{R}$-convex functions. Then the sum $\varphi_1(e_1 x)e_1 + \varphi_2(e_2 x)e_2$
is a $\mathbb{D}$-convex function.
\end{theorem}

\begin{proof}
Let $x ,y \in \mathbb{D}$. Then $e_1x, e_1y \in e_1\mathbb{D}$ and $e_2x, e_2y \in e_2\mathbb{D}.$ Further, let $0 \leq \lambda_1 \leq 1$ and $0 \leq \lambda_2 \leq 1$ be such that $0 \leq' \lambda \leq' 1$ where $\lambda= \lambda_1 e_1 + \lambda_2 e_2.$ Since $\varphi_1$ and $\varphi_2$ are $\mathbb{R}$- convex function, we have
\begin{equation}
\label{eq:N7}
\varphi_1(\lambda_1 e_1 x + (1-\lambda_1)e_1 y) \leq \lambda_1 \varphi_1(e_1 x)+(1-\lambda_1)\varphi_1(e_1y) 
\end{equation}
and 
\begin{equation}
\label{eq:N8}
\varphi_2(\lambda_2 e_2 x + (1-\lambda_2)e_2 y) \leq \lambda_2 \varphi_2(e_2 x)+(1-\lambda_2)\varphi_2(e_2y). 
\end{equation}
Now
\begin{eqnarray*}
\varphi_1(\lambda_1 e_1 x + (1-\lambda_1)e_1 y) & = & \varphi_1((\lambda_1e_1 + \lambda_2e_2) e_1 x + ((1-\lambda_1)e_1 + (1-\lambda_2)e_2)e_1 y)\\
& = & \varphi_1(\lambda e_1 x+(1-\lambda)e_1 y).\\
\end{eqnarray*}
So,
\begin{equation}
\label{eq:N9}
\varphi_1(\lambda_1 e_1 x + (1-\lambda_1)e_1 y)e_1= \varphi_1(\lambda e_1 x+(1-\lambda)e_1 y)e_1
\end{equation}
Similarly,
\begin{equation}
\label{eq:N10}
\varphi_2(\lambda_2 e_2 x + (1-\lambda_2)e_2 y)e_2= \varphi_2(\lambda e_2 x+(1-\lambda)e_2 y)e_2
\end{equation}
Further, it is easy to see that
\begin{equation}
\label{eq:N11}
\left[\lambda_1 \varphi_1(e_1 x)+(1-\lambda_1)\varphi_1(e_1y)\right]e_1=\left[\lambda \varphi_1(e_1 x)e_1 + (1-\lambda)\varphi_1(e_1y)e_1\right]
\end{equation}
and
\begin{equation}
\label{eq:N12}
\left[\lambda_2 \varphi_2(e_2 x)+(1-\lambda_2)\varphi_2(e_2y)\right]e_2=\left[\lambda \varphi_2(e_2 x)e_2 + (1-\lambda)\varphi_2(e_2y)e_2\right].
\end{equation}
Multiplying equation \eqref{eq:N7} and \eqref{eq:N8} by $e_1$ and $e_2$ and then adding we get
\begin{eqnarray*}
\varphi_1(\lambda_1 e_1 x + (1-\lambda_1)e_1 y)e_1 + \varphi_2(\lambda_2 e_2 x + (1-\lambda_2)e_2 y) e_2 & \leq' & \left[\lambda_1 \varphi_1(e_1 x)+(1-\lambda_1)\varphi_1(e_1y)\right] e_1\\
 & + & \left[\lambda_2 \varphi_2(e_2 x)+(1-\lambda_2)\varphi_2(e_2y)\right] e_2 \\
\end{eqnarray*} 
Using \eqref{eq:N9} ,~\eqref{eq:N10},~ \eqref{eq:N11} and \eqref{eq:N12} in the previous equation we finally get
\begin{eqnarray*}
\varphi_1(\lambda e_1 x+(1-\lambda)e_1 y)e_1 & + & \varphi_2(\lambda e_2 x+(1-\lambda)e_2 y)e_2\\
                        & \leq' &  \left[\lambda \varphi_1(e_1 x)e_1 + (1-\lambda)\varphi_1(e_1y)e_1\right]\\
												& + & \left[\lambda \varphi_2(e_2 x)e_2 + (1-\lambda)\varphi_2(e_2y)e_2\right] \\
\end{eqnarray*}
Simplifying we finally get,
$$\left[\varphi_1(e_1 . )e_1  +  \varphi_2(e_2 .)e_2\right](\lambda x+(1-\lambda)y) \leq' \lambda\left(\varphi_1(e_1.) + \varphi_1(e_2.)(x)\right)+(1-\lambda)\left(\varphi_1(e_1.) + \varphi_1(e_2.)(y)\right)$$
This proves that $\varphi_1e_1+\varphi_2e_2$ is also a $\mathbb{D}$- convex function.  
\end{proof}


\begin{theorem}(\cite{Mo-Ro-98}) If $\varphi(h)=f'(h)$ is hyperbolic derivable function and the Jordan arc $\gamma$ is contained in the domain of $\varphi$, then we have 
\begin{equation}
\label{eq:IN1}
\int_{\gamma} f(h)dh=\varphi(\gamma(b))-\varphi(\gamma(a))
\end{equation} 
\end{theorem}

\begin{theorem}
Let $\varphi_{\mathbb{D}}: \mathcal{D}_{\mathbb{D}}(h_0,r) \to \mathbb{D}$ be a $\mathbb{D}$-valued convex function, where $$\mathcal{D}_{\mathbb{D}}(h_0,r)=\{h \in \mathbb{D}:|h-h_0|_k \leq ' r \}$$ Then for every closed disk $\overline{\mathcal{D}_{\mathbb{D}}(y,\epsilon)} \subset \mathcal{D}_{\mathbb{D}}(h_0,r)$, $\varphi_{\mathbb{D}}$ has an integral representation as 
\begin{equation}
\varphi_{\mathbb{D}}(\gamma(x))=\varphi_{\mathbb{D}}(\gamma(h_0))+\int_{\Gamma} f(t)d\mu_{\mathbb{D}}(t)
\end{equation}
for every Jordan arc $\Gamma= Im\gamma$ with $\gamma:\[a,b\] \to \mathbb{D}$ and $\gamma'(t)$ is continuous except for a finite number of points and $f: \mathbb{D} \to \mathbb{D}$ is a $\mathbb{D}$-monotone non-decreasing and left continuous function. Also the derivative of $\varphi_{\mathbb{D}}$ along any line passing through each point in $\mathcal{D}_{\mathbb{D}}(h_0,r)$ exists and are equal except for a countable number of points. 
\end{theorem}

\begin{proof}
Suppose that $\varphi_{\mathbb{D}}$ is a $\mathbb{D}$-valued convex function and let $h_1,h_2,h_3 \in \mathcal{D}_{\mathbb{D}}(y,\epsilon)$ be such that $h_1 <' h_2 <' h_3$ under the partial ordering defined on $\mathbb{D}$ with $h_i-h_j$ for $i,j=1,2,3$ lies outside the null cone of $\mathbb{D}.$
Setting $\alpha=\frac{h_2-h_1}{h_3-h_1}$ and $\beta=\frac{h_3-h_2}{h_3-h_1}$. Then it is easy to see that $0<'\alpha, \beta <' 1,$ and $\alpha+\beta =1.$ Now with $x_1=h_3$ and $x_2=h_1$, we have $h_2=\alpha x_1+ \beta x_2$, so that 
\begin{eqnarray*}
\varphi_{\mathbb{D}}(h_2) & = & \varphi_{\mathbb{D}}(\alpha x_1 + \beta x_2) \\
                          & \leq' & \alpha \varphi_{\mathbb{D}}(x_1)+ \beta \varphi_{\mathbb{D}}(x_2) \\
													& = & \alpha \varphi_{\mathbb{D}}(h_3)+ \beta \varphi_{\mathbb{D}}(h_1).
\end{eqnarray*}
Substituting the values of $\alpha$ and $\beta$ in the previous equation, we get 
\begin{equation}
\label{eq:EC1}
(h_3-h_1)\varphi_{\mathbb{D}}(h_2)  \leq'  (h_2-h_1)\varphi_{\mathbb{D}}(h_3) + (h_3-h_2)\varphi_{\mathbb{D}}(h_1). 
\end{equation}
Solving equation \eqref{eq:EC1}, we get 
\begin{equation}
\label{eq:EC2}
\frac{\varphi_{\mathbb{D}}(h_2)-\varphi_{\mathbb{D}}(h_1)}{h_2-h_1} \leq' \frac{\varphi_{\mathbb{D}}(h_3)-\varphi_{\mathbb{D}}(h_1)}{h_3-h_1}. 
\end{equation}
Also we have 
\begin{equation}
\label{eq:EC3}
(h_3-h_1)\varphi_{\mathbb{D}}(h_2)  \leq  (h_3-h_2)\varphi_{\mathbb{D}}(h_1)+ ((h_3-h_1)-(h_3-h_2))\varphi_{\mathbb{D}}(h_3).  
\end{equation}
Solving equation \eqref{eq:EC3}, we get

\begin{equation}
\label{eq:EC4}
\frac{\varphi_{\mathbb{D}}(h_3)-\varphi_{\mathbb{D}}(h_1)}{h_3-h_1} \leq' \frac{\varphi_{\mathbb{D}}(h_3)-\varphi_{\mathbb{D}}(h_2)}{h_3-h_2}. 
\end{equation}
Combining equation \eqref{eq:EC2} and \eqref{eq:EC4}, we get
\begin{equation}
\label{eq:EC5}
\frac{\varphi_{\mathbb{D}}(h_2)-\varphi_{\mathbb{D}}(h_1)}{h_2-h_1} \leq' \frac{\varphi_{\mathbb{D}}(h_3)-\varphi_{\mathbb{D}}(h_1)}{h_3-h_1} \leq' \frac{\varphi_{\mathbb{D}}(h_3)-\varphi_{\mathbb{D}}(h_2)}{h_3-h_2}
\end{equation} 

This implies that the difference quotient is non-decreasing in $\mathcal{D}_{\mathbb{D}}(y,\epsilon)$. Now to prove the differentiability of $\varphi_{\mathbb{D}}$ at each point of $\mathcal{D}_{\mathbb{D}}(h_0,r) $. Let us consider an arbitrary point $x \in \mathcal{D}_{\mathbb{D}}(h_0,r)$ and let $\left\{L_i : i \in \Delta \right\}$ be the family of lines passing through $x$, where $\Delta$ is an indexing set.\\

To prove the differentiability of $\varphi_{\mathbb{D}}$  along any line say $L_k$ from two opposite direction of $x$ along $L_k.$ Let us first check the differentiability of $\varphi_{\mathbb{D}}$ along the right hand(positive direction) of $x$ along $L_k.$ Let us take $x_1,x_2$ be such that $x<' x_1 <' x_2$ and $x_i-x_j \notin \mathcal{NC}$ for $i,j=1,2$. Now identifying the points $x,x_1,x_2$ respectively with $h_1,h_2,h_3$ in \eqref{eq:EC5}, we get
\begin{equation}
\label{eq:EC6}
\frac{\varphi_{\mathbb{D}}(x_1)-\varphi_{\mathbb{D}}(x)}{x_1-x} \leq' \frac{\varphi_{\mathbb{D}}(x_2)-\varphi_{\mathbb{D}}(x)}{x_2-x} \leq' \frac{\varphi_{\mathbb{D}}(x_2)-\varphi_{\mathbb{D}}(x_1)}{x_2-x_1}
\end{equation}
This shows that~~$\frac{\varphi_{\mathbb{D}}(h)-\varphi_{\mathbb{D}}(x)}{h-x}$~~$\mathbb{D}$-decreasing as $h$ tends to $x$ from the positive direction along $L_k.$ So that the derivative
$$D^{L_k^+}_{\varphi_{\mathbb{D}}}(x)=\lim_{h \to x^+~~~along~~L_k} \frac{\varphi_{\mathbb{D}}(h)-\varphi_{\mathbb{D}}(x)}{h-x}$$ exists in extended hyperbolic plane $\overline{\mathbb{D}}$ and 
\begin{equation}
\label{eq:EC7}
D^{L_k^+}_{\varphi_{\mathbb{D}}}(x) \leq' \frac{\varphi_{\mathbb{D}}(h)-\varphi_{\mathbb{D}}(x)}{h-x}~~~~~\forall~~~h>' x~~~\mbox{and}~~~h-x \notin \mathcal{NC}
\end{equation}
where $\mathcal{NC}$ denotes the null cone of $\mathbb{D}.$ Next we take $u,h \in \mathcal{D}_{\mathbb{D}}(h_0,r)$ such that $u <' x <' h$ and difference of any two of this elements lies outside $\mathcal{NC}.$ Now identifying the points respectively with $h_1,h_2,h_3$ in \eqref{eq:EC5}, we obtain 

\begin{equation}
\frac{\varphi_{\mathbb{D}}(x)-\varphi_{\mathbb{D}}(u)}{x-u} \leq' \frac{\varphi_{\mathbb{D}}(h)-\varphi_{\mathbb{D}}(x)}{h-x}. 
\end{equation}
Thus for every $h>'x$, the quotient $\frac{\varphi_{\mathbb{D}}(h)-\varphi_{\mathbb{D}}(x)}{h-x}$ is $\mathbb{D}$-bounded below by some hyperbolic constant. Thus we have 
$$D^{L_k^+}_{\varphi_{\mathbb{D}}}(x)=\lim_{h \to x^+~~~along~~L_k} \frac{\varphi_{\mathbb{D}}(h)-\varphi_{\mathbb{D}}(x)}{h-x}$$ exists in hyperbolic plane $\mathbb{D}$, and hence $\varphi_{\mathbb{D}}$ is differentiable along the positive direction of the line $L_k$ passing through $x.$ Similarly $\varphi_{\mathbb{D}}$ is differentiable at $x$ along the left(negative) direction of $L_k.$ That is $D^{L_k^-}_{\varphi_{\mathbb{D}}}(x)$ exits and is $\mathbb{D}$-bounded above by some hyperbolic constant. Thus $D^{L_k^-}_{\varphi_{\mathbb{D}}}(x)$ exists in $\mathbb{D}.$. Next we shall show that $\varphi_{\mathbb{D}}$ satisfies the $\mathbb{D}$-Lipschitz condition.\\
\\
For this we see that the derivative of $\varphi_{\mathbb{D}}$ along any line in both the positive and negative direction exists at each point of $\mathcal{D}_{\mathbb{D}}(h_0,r)$ and for $h <' t$ we have 

\begin{equation}
\label{eq:EC8}
D^{L_i^+}_{\varphi_{\mathbb{D}}}(h) \leq' \frac{\varphi_{\mathbb{D}}(t)-\varphi_{\mathbb{D}}(h)}{t-h} \leq'  D^{L_i^-}_{\varphi_{\mathbb{D}}}(t),
\end{equation}
where $L_i$ is the line passing through $h,t$. Now for every $h_1,h_2$ with $h_1 < 'h_2$, we have 
\begin{equation}
\label{eq:EC9}
D^{L_i^+}_{\varphi_{\mathbb{D}}}(c) \leq' D^{L_i^+}_{\varphi_{\mathbb{D}}}(h_1) \leq' \frac{\varphi_{\mathbb{D}}(h_2)-\varphi_{\mathbb{D}}(h_1)}{h_2-h_1} \leq' D^{L_j^-}_{\varphi_{\mathbb{D}}}(h_2) \leq' D^{L_j^-}_{\varphi_{\mathbb{D}}}(d) 
\end{equation}
for the $\mathbb{D}$-smallest number $c$ and $\mathbb{D}$-largest number d such that $c \leq' h_1 \leq' h_2 \leq' d.$ Take 
$$M=max \left\{ |D^{L_i^+}_{\varphi_{\mathbb{D}}}(c)|_k, |D^{L_j^-}_{\varphi_{\mathbb{D}}}(d)|_k \right\}.$$
Then
\begin{equation}
\label{eq:EC10}
|\varphi_{\mathbb{D}}(h_2)-\varphi_{\mathbb{D}}(h_1)|_k \leq' M |h_2-h_1|_k
\end{equation}

Hence $\varphi_{\mathbb{D}}$ is $\mathbb{D}$-absolutely continuous and so by First Fundamental theorem for hyperbolic calculus, we have

\begin{equation}
\label{eq:EC11}
\varphi_{\mathbb{D}}(\gamma(x)) - \varphi_{\mathbb{D}}(\gamma(x))=\int_{\Gamma} \varphi_{\mathbb{D}}'(z)d\mu_{\mathbb{D}}
\end{equation}

Next finally we shall verify the properties of $\varphi_{\mathbb{D}}.$ Since $D^{L_i^+}$ and $D^{L_i^-}$ exists at each point of $\mathcal{D}_{\mathbb{D}}(h_0,r)$, we have for $h<'t$ \, from equation \eqref{eq:EC8} $D^{L_i^+}_{\varphi_{\mathbb{D}}}(h) \leq'  D^{L_i^-}_{\varphi_{\mathbb{D}}}(t)$. Thus we see that $D^{L_i^+\pm}(.)$ is $\mathbb{D}$-increasing function. So that the set of discontinuities of these function is almost countable and so $D^{L_i^+}(h)=D^{L_i^-}(h)$ at each continuous point of this function and this common values is $f=\varphi_{\mathbb{D}}'$. Thus we have from \eqref{eq:EC11} that 

\begin{equation}
\label{eq:EC11}
\varphi_{\mathbb{D}}(\gamma(x)) = \varphi_{\mathbb{D}}(\gamma(x))+\int_{\Gamma} f(t)d\mu_{\mathbb{D}}(t)
\end{equation}
\end{proof}

\begin{definition}($\mathbb{D}$-valued Young Function)
A convex function $\varphi_{\mathbb{D}}:\mathbb{D}^+ \to \overline{\mathbb{D}}^+$ is said to be an Orlicz function if it satisfies the following condition\\
(i) $\varphi_{\mathbb{D}}(0)=0.$\\
(ii) $\lim_{x \to \infty} \varphi_{\mathbb{D}}(x)= +\infty$ where we assume the convention that\\
$ +\infty=\alpha e+ \infty e_{2}=\infty e + \beta e_{2} = \infty e + \infty e_{2}$\\ and $\lim_{x \to \infty} \varphi_{\mathbb{D}}(x)$ means that the limit must exists along any line in the hyperbolic plane and must be equal.
 \end{definition}
Let us construct an example in this direction.
\begin{example}
Consider $\varphi_{\mathbb{D}}: \mathbb{D}^+ \to \mathbb{\overline{D}}^+$ defined by $\varphi_{\mathbb{D}}(x)=x^p$ for $p \geq 1.$
Then it is easy to check that $\varphi_{\mathbb{D}}$ is a $\mathbb{D}$-valued convex function. Further we see that
$$\varphi_{\mathbb{D}}(x)=x^p=\rho_1^p e_1+\rho_1^p e_2,~~\mbox{where}~~x=\rho_1 e_1 + \rho_2 e_2.$$
\end{example}

\begin{definition}
An Orlicz function $\varphi_{\mathbb{D}}:\mathbb{D}^+ \to \mathbb{D}^+$ is said to be a $\mathbb{D}$-valued N- function if 
\begin{itemize}
\item[(i)] $\varphi_{\mathbb{D}}$ is continuous.\\
\item[(ii)] $\varphi_{\mathbb{D}}(x)=0$ $\Leftrightarrow x=0.$\\
\item[(iii)] $\lim_{x \to 0} \frac{\varphi_{\mathbb{D}}(x)}{x}=0$~~$\forall$~~$x \in \mathbb{D}^+ \setminus \mathcal{NC}.$\\
\item[(iv)] $\lim_{x \to 0} \frac{\varphi_{\mathbb{D}}(x)}{x}=0$~~~$\forall~~x \in \mathbb{D}^+ \setminus \mathcal{NC}.$
\end{itemize}
\end{definition}

\begin{definition}
Let $\varphi_{\mathbb{D}}:\mathbb{D}^+ \to \mathbb{\overline{D}}^+$ be a $\mathbb{D}$-Young function. We define another Young function $\psi_{\mathbb{D}}:\mathbb{D}^+ \to \mathbb{\overline{D}}^+$ as follows
$$\psi_{\mathbb{D}}(y)=\sup\left\{xy-\varphi_{\mathbb{D}}(x): x \in \mathbb{D}^+\right\}.$$
We call $\psi_{\mathbb{D}}$ a $\mathbb{D}$-complementary function corresponding to $\varphi_{\mathbb{D}}$.
\end{definition}
\begin{remark}
\label{Re:1}
Let $x=\rho_1 e_1+\rho_2 e_2$ and $y=\gamma_1 e_1 + \gamma_2 e_2$ be two elements of $\mathbb{D}^+$. Then we see that
$xy=\rho_1\gamma_1 e_1+\rho_2\gamma_2 e_2$ and since $\varphi_{\mathbb{D}}(x)= \varphi_{\mathbb{D}_1}(e_1x)e_1 + \varphi_{\mathbb{D}_2}(e_2x)e_2$, we have 
\begin{eqnarray*}
xy-\varphi_{\mathbb{D}}(x) & = & \left(\rho_1\gamma_1 e_1+\rho_2\gamma_2 e_2\right)-\left(\varphi_{\mathbb{D}_1}(e_1x)e_1 + \varphi_{\mathbb{D}_2}(e_2 x)e_2\right)\\
&=& \left(\rho_1\gamma_1 - \varphi_{\mathbb{D}_1}(e_1x)\right)e_1 + \left(\rho_2\gamma_2- \varphi_{\mathbb{D}_2}(e_2x)\right)e_2\\
\end{eqnarray*}
Now taking supermum over all $x \in \mathbb{D}^+$, we obtain
\begin{eqnarray*}
\psi_{\mathbb{D}}(y)=\sup\left\{ xy-\varphi_{\mathbb{D}}(x): x \in \mathbb{D}^+\right\}= \sup_{\rho_1 \geq 0}\left\{ \rho_1\gamma_1 - \varphi_{\mathbb{D}_1}(e_1x)\right\}e_1+ \sup_{\rho_2 \geq 0}\left\{\rho_2\gamma_2- \varphi_{\mathbb{D}_2}(e_2x)\right\}e_2
\end{eqnarray*}
\end{remark}

\begin{theorem}
Let $\varphi_{\mathbb{D}}$ be a Young function and $\psi_{\mathbb{D}}$ be the $\mathbb{D}$-complementary pair of $\varphi_{\mathbb{D}}$. Then the pair $\left(\varphi_{\mathbb{D}}, \psi_{\mathbb{D}}\right)$ satisfies the Young's inequality 
$$xy \leq ' \varphi_{\mathbb{D}}(x)+ \psi_{\mathbb{D}}(y).$$
\end{theorem}
\begin{proof}
Since $\psi_{\mathbb{D}}$ is the $\mathbb{D}$-complementary pair of $\varphi_{\mathbb{D}}$ and so we have from Remark \ref{Re:1} 
\begin{eqnarray*}
\psi_{\mathbb{D}}(y)= \sup_{\rho_1 \geq 0}\left\{ \rho_1\gamma_1 - \varphi_{\mathbb{D}_1}(e_1x)\right\}e_1+ \sup_{\rho_2 \geq 0}\left\{\rho_2\gamma_2- \varphi_{\mathbb{D}_2}(e_2x)\right\}e_2.
\end{eqnarray*}
Further since $\psi_{\mathbb{D}}:\mathbb{D}^+ \to \mathbb{D}^+$ is a Young function, we have
$$\psi_{\mathbb{D}}(y)=\psi_{\mathbb{D}_1}(e_1y)e_1+ \psi_{\mathbb{D}_2}(e_2y)e_2.$$
So,
$\psi_{\mathbb{D}_1}(e_1y)= \sup_{\rho_1 \geq 0}\left\{ \rho_1\gamma_1 - \varphi_{\mathbb{D}_1}(e_1x)\right\}$ and 
$\psi_{\mathbb{D}_2}(e_2y)=\sup_{\rho_2 \geq 0}\left\{\rho_2\gamma_2- \varphi_{\mathbb{D}_2}(e_2x)\right\}.$
Since $(\varphi_{\mathbb{D}_1},\psi_{\mathbb{D}_1})$ and $(\varphi_{\mathbb{D}_2},\psi_{\mathbb{D}_2})$ are two $\mathbb{R}$-complementary pair of Young functions, we have
\begin{equation}
\label{eq:2e}
\rho_1\gamma_1 - \varphi_{\mathbb{D}_1}(e_1x) \leq \psi_{\mathbb{D}_1}(e_1y)~~ \mbox{and}~~ \rho_2\gamma_2- \varphi_{\mathbb{D}_2}(e_2x) \leq \psi_{\mathbb{D}_2}(e_2y).
\end{equation}
Using equation \eqref{eq:2e}, we get 
$$\rho_1\gamma_1e_1 + \rho_2\gamma_2e_2 \leq ' \left(\varphi_{\mathbb{D}_1}(e_1x)e_1 + \varphi_{\mathbb{D}_2}(e_2x) \right) + \left(\psi_{\mathbb{D}_1}(e_1y)e_1 + \psi_{\mathbb{D}_2}(e_2y)e_2 \right).$$
Thus,
$$xy \leq' \varphi_{\mathbb{D}}(x)+ \psi_{\mathbb{D}}(y).$$
This proves the Young inequality for a $\mathbb{D}$-complementary pair of Young functions. 
\end{proof}

\begin{definition}
A $\mathbb{D}$-valued Young function $\varphi_{\mathbb{D}}:\mathbb{D}^+ \to \mathbb{D}^+$ is said to satisfy the $\Delta_{\mathbb{D}}^2$-condition denoted by $\varphi_{\mathbb{D}} \in \Delta_{\mathbb{D}}^2$ if there exist some hyperbolic constants $K \geq' 0$ and $x_o$(depending upon K) such that
$$\varphi_{\mathbb{D}}((2e_1+2e_2)x) \leq' K\varphi_{\mathbb{D}}(x)~~~~~~~\forall~~~0 \leq' x \leq' x_o.$$
\end{definition}
\begin{corollary}
Let $\varphi_{\mathbb{D}}$ be a $\mathbb{D}$-valued Young function, Then $\varphi_{\mathbb{D}} \in \Delta_{\mathbb{D}}^2$ if and only if $\varphi_{\mathbb{D}_1} \in \Delta^2$ and $\varphi_{\mathbb{D}_1} \in \Delta^2.$
\end{corollary}

\section{ \bf $\mathbb{D}$-Metric Spaces.}
In this section we shall define $\mathbb{D}$-metric spaces and prove the completeness of a $\mathbb{D}$-metric space.
\begin{definition}
Let $\Omega_{\mathbb{D}}$ be a non-empty set. Then a hyperbolic valued function
$$d_{\mathbb{D}}:\Omega_{\mathbb{D}} \times \Omega_{\mathbb{D}} \to \mathbb{D}$$
is called a $\mathbb{D}$-metric on $\Omega_{\mathbb{D}}$, if for every $x,y,z \in \Omega_{\mathbb{D}}$, the following holds.\\
\begin{itemize}
\item[(i)] $d_{\mathbb{D}}(x,y) \geq' 0$ and $d_{\mathbb{D}}(x,y)=0$ if and only if $x=y$.\\
\item[(ii)] $d_{\mathbb{D}}(x,y)=d_{\mathbb{D}}(y,x).$\\
\item[(iii)] $d_{\mathbb{D}}(x,z) \leq' d_{\mathbb{D}}(x,y)+ d_{\mathbb{D}}(y,z).$  
\end{itemize}
The pair $\left(\Omega_{\mathbb{D}}, d_{\mathbb{D}}\right)$ is called a hyperbolic metric space and we say that $d_{\mathbb{D}}$ is a $\mathbb{D}$-metric on $\Omega_{\mathbb{D}}$ and $\left(\Omega_{\mathbb{D}}, d_{\mathbb{D}}\right)$ is a $\mathbb{D}$ -metric space.
\end{definition}

Now since $d_{\mathbb{D}}:\Omega_{\mathbb{D}} \times \Omega_{\mathbb{D}} \to \mathbb{D}$ is a $\mathbb{D}$-valued metric, we can   decompose  $d_\mathbb{D}$ as follows
$$d_\mathbb{D}(x,y)= d_{\mathbb{D}}^1(x,y)e_1+d_{\mathbb{D}}^2(x,y)e_2.$$
where $d_{\mathbb{D}}^1$ and $d_{\mathbb{D}}^2$ are real valued metric defined on $\Omega_{\mathbb{D}} \times \Omega_{\mathbb{D}}$. We say that $d_{\mathbb{D}}^1$ and $d_{\mathbb{D}}^2$ are $\mathbb{R}$-metrics.

\begin{definition}
Let $\Omega_{\mathbb{D}}$ be a non empty set and $d_{\mathbb{D}}$ be a $\mathbb{D}$-metric defined on $\Omega_{\mathbb{D}}$, so that $\left(\Omega_{\mathbb{D}},d_{\mathbb{D}}\right)$ is a $\mathbb{D}$-metric space. Then\\
\begin{itemize}
\item[(i)] A sequence $\left\{x_n\right\}$ of points in $\Omega_{\mathbb{D}}$ is said to converge to a point 
$x_o \in \Omega_{\mathbb{D}}$ with respect to the $\mathbb{D}$-metric $d_{\mathbb{D}}$, if for every $\mathbb{D}$-small hyperbolic number $\epsilon >' 0$, there exists $m \in \mathbb{N}$ such that
$$d_{\mathbb{D}}(x_n,x_o) <' \epsilon~~~~~~~~~~\forall~~n>m.$$
\item[(ii)] A sequence $\left\{x_n\right\}$ of points in $\Omega_{\mathbb{D}}$ is said to be a hyperbolic Cauchy sequence if for each $\mathbb{D}$-small hyperbolic number $\epsilon >' 0$, there exists a positive integer $N$ such that
$$d_{\mathbb{D}}(x_n,x_m) <' \epsilon~~~~~~~~~~\forall~~n,m>N.$$
\item[(iii)] Let $\left(\Omega_{\mathbb{D}}, d_{\mathbb{D}}\right)$ be a $\mathbb{D}$-metric space. Then $\left(\Omega_{\mathbb{D}}, d_{\mathbb{D}}\right)$ is said to be a complete $\mathbb{D}$-metric space if every hyperbolic Cauchy sequence in $\Omega_{\mathbb{D}}$ converges in $\Omega_{\mathbb{D}}.$
\end{itemize}
\end{definition}

\begin{theorem}
Let $\Omega_{\mathbb{D}}$ be a $\mathbb{D}$-module. Then
$$d_\mathbb{D}= d_{\mathbb{D}}^1e_1+d_{\mathbb{D}}^2e_2.$$
is a $\mathbb{D}$-metric on $\Omega_{\mathbb{D}}$ if and only if $d_{\mathbb{D}}^1$ and $d_{\mathbb{D}}^2$ are $\mathbb{R}$-metric on $\Omega_{\mathbb{D}}$.
\end{theorem}

\begin{theorem}
Let $\left(\Omega_{\mathbb{D}},d_{\mathbb{D}}\right)$ be a complete  $\mathbb{D}$-metric spaces. Then $\left(\Omega_{\mathbb{D}},d_{\mathbb{D}}^1\right)$ and $\left(\Omega_{\mathbb{D}},d_{\mathbb{D}}^2\right)$ are complete $\mathbb{R}$-metric spaces if $d_{\mathbb{D}}^1$ and $d_{\mathbb{D}}^2$ are equivalent metrics on $\Omega_{\mathbb{D}}$, where $d_{\mathbb{D}}^1$ and $d_{\mathbb{D}}^2$ is such that 
$$d_{\mathbb{D}}=d_{\mathbb{D}}^1e_1+d_{\mathbb{D}}^2 e_2.$$
\end{theorem}

\begin{proof}
Suppose that $\left(\Omega_{\mathbb{D}},d_{\mathbb{D}}\right)$ is a complete  $\mathbb{D}$-metric space. To show that $\left(\Omega_{\mathbb{D}},d_{\mathbb{D}}^1\right)$ and $\left(\Omega_{\mathbb{D}},d_{\mathbb{D}}^2\right)$ are complete $\mathbb{R}$-metric spaces. Consider a Cauchy sequence in any one of the $\mathbb{R}$-metric space, as $d_{\mathbb{D}}^1$ and   $d_{\mathbb{D}}^2$ are equivalent metric on $\Omega_{\mathbb{D}}$ and any Cauchy sequence in $\left(\Omega_{\mathbb{D}},d_{\mathbb{D}}^1\right)$ will also be a Cauchy sequence in $\left(\Omega_{\mathbb{D}},d_{\mathbb{D}}^2\right)$ and conversely.  So let  us suppose that $\{x_n\}$ be a Cauchy sequence in $\left(\Omega_{\mathbb{D}},d_{\mathbb{D}}^1\right)$ and $\left(\Omega_{\mathbb{D}},d_{\mathbb{D}}^2\right).$ This means that for every $\epsilon_1 >0,~\exists~ N_1~ \in \mathbb{N}$ and for every $\epsilon_2 >0,~\exists~ N_2 \in \mathbb{N}$ such that 
\begin{equation}
\label{eq:M1}
d_{\mathbb{D}}^1(x_n,x_m) < \epsilon_1~~~\forall~~~n,m>N_1~~~\mbox{and}~~ d_{\mathbb{D}}^2(x_n,x_m) < \epsilon_2~~~\forall~~~n,m>N_2.
\end{equation}
Using equation \eqref{eq:M1}, we get 
\begin{equation}
d_{\mathbb{D}}(x_n,x_m) <' \epsilon~~~\forall~~~n,m>N
\end{equation}
where $\epsilon=\epsilon_1e_1+\epsilon_2e_2$ and $N=max\{N_1,N_2\}.$
Thus $\{x_n\}$ becomes a Cauchy sequence in $\left(\Omega_{\mathbb{D}},d_{\mathbb{D}}\right)$. By completeness of $\mathbb{D}$-metric spaces $\left(\Omega_{\mathbb{D}},d_{\mathbb{D}}\right)$, there exists some $y \in \Omega_{\mathbb{D}}$,  such that for every $\mathbb{D}$-small hyperbolic number  $\delta>'0,~\exists~ M \in \mathbb{N}$ such that 
\begin{equation}
d_{\mathbb{D}}(x_n,y) <' \delta ~~~\forall~~~n >M.
\end{equation}
This implies that 
$$d_{\mathbb{D}}^1(x_n,y)e_1+d_{\mathbb{D}}^2(x_n,y)e_2 <' \delta_1e_1+\delta_2e_2 ~~~\forall~~~n >M.$$
So, $d_{\mathbb{D}}^1(x_n,y) < \delta_1$, and $d_{\mathbb{D}}^2(x_n,y) < \delta_2$~$\forall~ n > M$  Thus, the Cauchy sequence $\{x_n\}$ in $\left(\Omega_{\mathbb{D}},d_{\mathbb{D}}^1\right)$ and $\left(\Omega_{\mathbb{D}},d_{\mathbb{D}}^2\right)$ converges to the point $y \in \Omega_{\mathbb{D}}$. This proves the completeness of the $\mathbb{R}$-metric spaces.
\end{proof}

\begin{theorem}
Let $\Omega_{\mathbb{D}}$ be any arbitrary non empty set. Let $d_{\mathbb{D}}^1$ and $d_{\mathbb{D}}^2$ be two $\mathbb{R}$-metrics defined on $\Omega_{\mathbb{D}}$. Suppose that $\left(\Omega_{\mathbb{D}},d_{\mathbb{D}}^1\right)$ , $\left(\Omega_{\mathbb{D}},d_{\mathbb{D}}^2\right)$ be complete $\mathbb{R}$-metric spaces. Then $\left(\Omega_{\mathbb{D}},d_{\mathbb{D}}\right)$ is a complete $\mathbb{D}$-metric spaces if $d_{\mathbb{D}}^1$ and $d_{\mathbb{D}}^2$ are two equivalent metric on $\Omega_{\mathbb{D}},$ where
\begin{equation}
d_{\mathbb{D}}= d_{\mathbb{D}}^1 e_1 +d_{\mathbb{D}}^2 e_2
\end{equation}
\end{theorem}

\begin{proof}
Suppose that $\left\{x_n\right\}$ be a Cauchy sequence in the $\mathbb{D}$-metric space $\left(\Omega_{\mathbb{D}},d_{\mathbb{D}}\right).$ Then for every $\epsilon >'0$ ($\mathbb{D}$-small hyperbolic number), there exists $N \in \mathbb{N}$ such that
\begin{equation}
d_{\mathbb{D}}(x_n,x_m) <' \epsilon~~\forall~~n,m>N.
\end{equation}
This implies that 
\begin{equation}
d_{\mathbb{D}}^1(x_n,x_m) e_1 + d_{\mathbb{D}}^2(x_n,x_m)e_2 <' \epsilon,~~\forall~~n,m>N.
\end{equation}
So that 
\begin{equation}
d_{\mathbb{D}}^1(x_n,x_m) < \epsilon_1~~~\mbox{and}~~~d_{\mathbb{D}}^2(x_n,x_m)e_2 < \epsilon_2,~~\forall~~n,m>N.
\end{equation}
where $\epsilon= \epsilon_1e_1+\epsilon_2e_2 >' 0.$
This show that $\left\{x_n\right\}$ is also a Cauchy sequence in the $\mathbb{R}$-metric spaces $\left(\Omega_{\mathbb{D}},d_{\mathbb{D}}^1\right)$ and $\left(\Omega_{\mathbb{D}} ,d_{\mathbb{D}}^2\right)$. Now since $\left(\Omega_{\mathbb{D}} ,d_{\mathbb{D}}^1\right)$ and $\left(\Omega_{\mathbb{D}} ,d_{\mathbb{D}}^2\right)$ are complete $\mathbb{R}$-metric spaces implies that $\{x_n\}$ converges to some point say $x$ in $\Omega_{\mathbb{D}}$ with respect to the $\mathbb{R}$-metric $d_{\mathbb{D}}^1$, and converges to some point say $y$ in $\Omega_{\mathbb{D}}$ with respect to the $\mathbb{R}$-metric $d_{\mathbb{D}}^2$. But since $d_{\mathbb{D}}^1$ and $d_{\mathbb{D}}^2$ are equivalent $\mathbb{R}$-metrics on $\Omega_{\mathbb{D}}$ and so this guarantees us with $x=y.$ This means that for every $\delta_1>0$ and $\delta_2>0$, there exists $N_1 , N_2 \in \mathbb{N}$ respectively such that 

\begin{equation}
\label{eq:C1}  
d_{\mathbb{D}}^1(x_n,x) < \delta_1~~\forall~~n>N_1~~\mbox{and}~~d_{\mathbb{D}}^2(x_n,x) < \delta_1~~\forall~~n>N_2
\end{equation}
Let $N=max\left\{N_1,N_2\right\}$ and from equation \eqref{eq:C1}, we get

\begin{equation}
d_{\mathbb{D}}^1(x_n,x)e_1+d_{\mathbb{D}}^2(x_n,x)e_2 <' \delta_1 e_1+ \delta_2 e_2~~~\forall~~n>N
\end{equation}
This proves that $\{x_n\}$ converges to the point $ x \in \mathbb{D}$ with respect to the metric $d_{\mathbb{D}}^1e_1+ d_{\mathbb{D}}^2e_2$  This proves the completeness of $\mathbb{D}$-metric space $\left(\Omega_{\mathbb{D}},d_{\mathbb{D}}^1e_1+ d_{\mathbb{D}}^2e_2\right).$ 

\end{proof}

\section{\bf Bicomplex Valued Function Spaces}

Let $\left(\Omega_{\mathbb{D}},\Sigma_{\mathbb{D}},\mu_{\mathbb{D}}\right)$ be a $\sigma$-finite $\mathbb{D}$- measure space. Consider the collection $\mathcal{M}(\left(\Omega_{\mathbb{D}},\Sigma_{\mathbb{D}},\mu_{\mathbb{BC}}\right))$ of all bicomplex  valued $\Sigma_{\mathbb{D}}$-measurable functions defined on $\Omega_{\mathbb{D}}$.
$$i.e.,~~~~\mathcal{M}\left(\Omega_{\mathbb{D}}\right)=\left\{f:\Omega_{\mathbb{D}} \to \mathbb{BC}:~ f~~\mbox{is}~~\Sigma_{\mathbb{D}}-\mbox{measurable function.}~\right\}$$
Thus if $f \in \mathcal{M}\left(\Omega_{\mathbb{D}}\right)$, then it is clear that we can decompose $f$ as
$f=f_1e_1+f_1e_2$, where $f_1=e_1 f$ and $f_2=e_2 f$ are complex valued $\Sigma_{\mathbb{D}}$-measurable functions on the measure spaces $\left(\Omega_{\mathbb{D}},\Sigma_{\mathbb{D}},\mu_{\mathbb{D}}^1\right)$ and $\left(\Omega_{\mathbb{D}},\Sigma_{\mathbb{D}},\mu_{\mathbb{D}}^2\right)$, respectively.\\ Under the operation of addition and scalar multiplication defined as
$$f+g=(f_1+g_1)e_1+(f_2+g_2)e_2~~~~\mbox{and}~~\alpha f=(\alpha_1f_1)e_1+(\alpha_2f_2)e_2.$$
for $f,g \in \mathcal{M}\left(\Omega_{\mathbb{D}}\right)$ and $\alpha \in \mathbb{BC}$, the collection $\left(\mathcal{M}\left(\Omega_{\mathbb{D}}\right),+,.\right)$ becomes a $\mathbb{BC}$ -module.\\
\\
Given a $\mathbb{D}$-measure space $\left(\Omega_{\mathbb{D}},\Sigma_{\mathbb{D}},\mu_{\mathbb{D}}\right).$ Two bicomplex valued $\Sigma_{\mathbb{D}}$-measurable functions $f$ and $g$ defined on a set $E \in \Sigma_{\mathbb{D}}$ are said to be $\mu_{\mathbb{D}}$-equal almost everywhere and we write
$$f=g~~~~\mbox{a.e on E,}~~$$
if there exists a set $N\subset E$ with $\mu_{\mathbb{D}}(N)=0$ such that $f(t)=g(t)$ for every $t \in E \setminus N.$
Further, since $f=f_1e_1+f_2e_2$ and $g=g_1e_1+g_2e_2$ implies that $f=g$ $\mu_{\mathbb{D}}$-a.e if and only if $f_1=g_1$ $\mu_{\mathbb{D}}^1$-a.e and  $f_2=g_2$ $\mu_{\mathbb{D}}^2$-a.e. Clearly the relation {\em equality almost everywhere} of two $\mathbb{BC}$-valued $\Sigma_{\mathbb{D}}$-measurable function is an equivalence relation and it partitions the collection $\mathcal{M}(\left(\Omega_{\mathbb{D}},\Sigma_{\mathbb{D}},\mu_{\mathbb{D}}\right))$ into disjoint equivalence classes. We denote by $\mathcal{M}^o(\left(\Omega_{\mathbb{D}},\Sigma_{\mathbb{D}},\mu_{\mathbb{D}}\right))$, the collection of all equivalence class of $\mathbb{BC}$-valued $\Sigma_{\mathbb{D}}$-measurable function which are identified $\mu_{\mathbb{D}}$-almost everywhere.

\begin{definition}({$\mathbb{D}$-Essentially Bounded functions}) Let $\left(\Omega_{\mathbb{D}},\Sigma_{\mathbb{D}},\mu_{\mathbb{D}}\right)$ be a $\mathbb{D}$-measure space and $f$ be an bicomplex valued measurable function defined on $\Omega_{\mathbb{D}}$. The essential $\mathbb{D}$-supremum of $f$ on $\Omega_{\mathbb{D}}$ denoted by $\|f\|_{\mathbb{D}}^{\infty}$ is defined as the infimum of all essential $\mathbb{D}$-bounds of $f$ on $\Omega_{\mathbb{D}}.$
$$i.e., \|f\|_{\mathbb{D}}^{\infty}=\inf\{\mathcal{K} \geq' 0 : \mu_{\mathbb{D}}(E)=0\}$$
where $E=\{t \in \Omega_{\mathbb{D}}: |f(t)|_k >' \mathcal{K} \}.$\\
Note that we take the infimum of empty set to be equal to $+\infty$ or $\infty e_1+be_2$ or $ae_1+\infty e_2.$
\end{definition}
\begin{remark}
Let us denote $L^{\infty,\mathbb{D}}(\Omega_{\mathbb{D}})$ to be the $\mathbb{BC}$- module(under the operation of addition and scalar multiplication of functions) of all equivalence class of bicomplex valued $\mathbb{D}$-essentially bounded, $\Sigma_{\mathbb{D}}$-measurable function that are identified $\mu_{\mathbb{D}}$-equal almost everywhere. Further every $f \in L^{\infty,\mathbb{D}}(\Omega_{\mathbb{D}})$  can be decomposed into
$f=f_1e_1+f_2e_2$ with $$\|f\|_{\mathbb{D}}^{\infty}=\|f_1e_1+f_2e_2\|_{\mathbb{D}}^{\infty}=\|f_1\|_{e_1\mathbb{D}}^{\infty} e_1+\|f_2\|_{e_2\mathbb{D}}^{\infty} e_2.$$
where $\|f_1\|_{e_1\mathbb{D}}^{\infty}$ and $\|f_2\|_{e_2\mathbb{D}}^{\infty}$ denote the essential bounds of complex valued $\Sigma_{\mathbb{D}}$-measurable function $f_1$ and $f_2.$ 
\end{remark} 
\begin{definition}{($\mathbb{D}$-Lebesgue integrable functions)}
Let $\left(\Omega_{\mathbb{D}},\Sigma_{\mathbb{D}},\mu_{\mathbb{D}}\right)$ be a $\mathbb{D}$-finite measure space and $f$ be a $\mathbb{D}$-bounded bicomplex valued $\Sigma_{\mathbb{D}}$-measurable function defined on $\Omega_{\mathbb{D}}$. Then the $\mathbb{D}$-Lebesgue integral of $f$ is defined as
\begin{eqnarray*}
\int_{\Omega_{\mathbb{D}}}|f|_kd\mu_{\mathbb{D}} & = & \int_{\Omega_{\mathbb{D}}}(|f_1|e_1+|f_2|e_2)d(\mu_{\mathbb{D}}^1 e_1+\mu_{\mathbb{D}}^2 e_2)\\
\end{eqnarray*}
and we  have the idempotent decomposition of this integral as
\begin{eqnarray*}
\int_{\Omega_{\mathbb{D}}}fd\mu_{\mathbb{D}} &=& \left(\int_{\Omega_{\mathbb{D}}}|f_1|d\mu_{\mathbb{D}}^1\right)e_1-\left(\int_{\Omega_{\mathbb{D}}}|f_2|d\mu_{\mathbb{D}}^1\right)e_2.\\
\end{eqnarray*}
where $f_1=e_1f$ and $f_2=e_2f$ are complex valued $\Sigma_{\mathbb{D}}$-measurable function.
\end{definition}

\begin{definition}
Let $\left(\Omega_{\mathbb{D}},\Sigma_{\mathbb{D}},\mu_{\mathbb{D}}\right)$ be a $\mathbb{D}$-finite measure space and let $\mathcal{A}$ be an arbitrary family of all $\mathbb{BC}$-valued $\Sigma_{\mathbb{D}}$-measurable function on $\Omega_{\mathbb{D}}$.
$$i.e.,~~~~~\mathcal{A}=\left\{f_n:\Omega_{\mathbb{D}} \to \mathbb{BC}, n \in \mathbb{N}\right\}.$$ Then the collection $\mathcal{A}$ is said to be $\mathbb{D}$-uniformly integrable, if the following condition holds:\\

\begin{itemize}
\item[(i)] $\sup_{n} \int_{\Omega_{\mathbb{D}}}|f_n|_kd\mu_{\mathbb{D}}= \sup_{n}\left(\int_{\Omega_{\mathbb{D}}}|f_n^1|d\mu_{\mathbb{D}}^1\right)e_1- \sup_{n}\left(\int_{\Omega_{\mathbb{D}}}|f_n^2|d\mu_{\mathbb{D}}^2\right)e_2 = C <' \infty.$\\
\item[(ii)] $\lim_{\mu_{\mathbb{D}}(A) \to 0} \int_{A}|f_n|_kd\mu_{\mathbb{D}}= \lim_{\mu_{\mathbb{D}}^1(A) \to 0}\left(\int_{\Omega_{\mathbb{D}}}|f_n^1|d\mu_{\mathbb{D}}^1\right)e_1- \lim_{\mu_{\mathbb{D}}^2(A) \to 0}\left(\int_{\Omega_{\mathbb{D}}}|f_n^2|d\mu_{\mathbb{D}}^2\right)e_2=0. $
\end{itemize}
Here $f_n^1=e_1f_n$ and $f_n^2=e_2f_n$ are complex valued $\Sigma_{\mathbb{D}}$-measurable functions.
\end{definition}
\begin{remark}
From above definition, we see that the collection $\mathcal{A}$ is $\mathbb{D}$-uniformly integrable if for every $f_n=f_n^1e_1+f_n^2e_2$ in $\mathcal{A}$, we must have\\ 
\begin{itemize}
\item[(i)] $\sup_{n}\left(\int_{\Omega_{\mathbb{D}}}|f_n^1|d\mu_{\mathbb{D}}^1\right) < \infty$ and $\sup_{n}\left(\int_{\Omega_{\mathbb{D}}}|f_n^2|d\mu_{\mathbb{D}}^2\right) < \infty.$ 
\item[(ii)] $\lim_{\mu_{\mathbb{D}}^1(A) \to 0}\left(\int_{\Omega_{\mathbb{D}}}|f_n^1|d\mu_{\mathbb{D}}^1\right)=0$
and $\lim_{\mu_{\mathbb{D}}^2(A) \to 0}\left(\int_{\Omega_{\mathbb{D}}}|f_n^2|d\mu_{\mathbb{D}}^2\right)=0.$ 
\end{itemize}
\end{remark}

Next we shall prove some equivalent conditions for $\mathcal{A}$ to be $\mathbb{D}$-uniformly integrable.
\begin{theorem} 
\label{th:38}
Consider the family $\mathcal{A}$ of all $\mathbb{BC}$-valued $\Sigma_{\mathbb{D}}$-measurable functions on a $\mathbb{D}$-finite measure space $\left(\Omega_{\mathbb{D}},\Sigma_{\mathbb{D}},\mu_{\mathbb{D}}\right).$ Then the following are equivalent:\\ 
\begin{itemize}
\item[(i)] $\mathcal{A}$ is $\mathbb{D}$-uniformly integrable.\\
\item[(ii)] $\lim_{\lambda \to 0} \int_{\left\{ |f_n|_k >'\lambda \right\}} |f_n|_k d\mu_{\mathbb{D}}=0$ $\mathbb{D}$-uniformly in $n \in \mathbb{N}.$
\item[(iii)] There always exists a $\mathbb{D}$-convex function $\varphi_{\mathbb{D}} :\mathbb{D}^+ \to  \mathbb{D}^+$ with $\varphi_{\mathbb{D}}(0)=0$ and  $\frac{\varphi_{\mathbb{D}}(x)}{x} \to \infty$ as $x \to \infty$ for every $x \in \mathbb{D} \setminus \mathcal{NC}$, in terms of which 
$$\sup_{n} \int_{\Omega_{\mathbb{D}}}\varphi_{\mathbb{D}}(|f_n|_k)d\mu_{\mathbb{D}} <' \infty.$$
\end{itemize}
\end{theorem}
\begin{proof}
Since the proof is well known for uniform integrability of real valued measurable function (see \cite{Rao-Ren91}) and using the idempotent decomposition for each $f_n \in \mathcal{A}$ as $|f_n|_k=|f_{n}^1|e_1+|f_{n}^2|e_2,$ we see that each of the above equivalent conditions are true for $\mathbb{BC}$-valued $\Sigma_{\mathbb{D}}$-measurable function if and only if it is true for each of the complex valued $\Sigma_{\mathbb{D}}$-measurable functions $f_n^1$ and $f_n^2.$ 
\end{proof}

\begin{definition}($\mathbb{BC}$-Orlicz Function Spaces.)
Let $\left(\Omega_{\mathbb{D}},\Sigma_{\mathbb{D}},\mu_{\mathbb{D}}\right)$ be a $\mathbb{D}$-finite measure space. Let $\varphi_{\mathbb{D}}: \mathbb{D}^+ \to \mathbb{D}^+$ be an Young function. Consider the collection  denoted by $\mathcal{F}_{\mathbb{D}}^{\varphi_{\mathbb{D}}}$ of all $\mathbb{BC}$-valued $\Sigma_{\mathbb{D}}$-measurable functions $f: \Omega_{\mathbb{D}} \to \mathbb{BC}$ satisfying the condition that $\int_{\mathbb{D}}\varphi_{\mathbb{D}}(|f|_k) <' \infty.$ 
$$i.e.,~~~~\mathcal{F}_{\mathbb{D}}^{\varphi_{\mathbb{D}}}=\left\{f:\Omega_{\mathbb{D}} \to \mathbb{BC}: f \mbox{ is $\Sigma_{\mathbb{D}}$-measurable and} \int_{\mathbb{D}}\varphi_{\mathbb{D}}(|f|_k)d\mu_{\mathbb{D}} <' \infty \right\}.$$

Further the idempotent decomposition defined for $\varphi_{\mathbb{D}}$ in \eqref{eq:1} gives  
\begin{equation}
\label{eq:23}
\int_{\Omega_{\mathbb{D}}}\varphi_{\mathbb{D}}(|f|_k)d\mu_{\mathbb{D}} = \left(\int_{\Omega_{\mathbb{D}}}\varphi_{\mathbb{D}_1}(|e_1f|)d\mu_{\mathbb{D}}^1\right)e_1- \left(\int_{\Omega_{\mathbb{D}}}\varphi_{\mathbb{D}_1}(|e_2f|)d\mu_{\mathbb{D}}^2\right)e_2.\\
\end{equation}
\end{definition}

\begin{remark} It is easy to see that the collection $\mathcal{F}_{\mathbb{D}}^{\varphi_{\mathbb{D}}}\left(\Omega_{\mathbb{D}},\Sigma_{\mathbb{D}},\mu_{\mathbb{D}}\right)$  is a $\mathbb{BC}$-linear module, and by using equation  \eqref{eq:23}, we have the idempotent decomposition of $\mathcal{F}_{\mathbb{D}}^{\varphi_{\mathbb{D}}}\left(\Omega_{\mathbb{D}},\Sigma_{\mathbb{D}},\mu_{\mathbb{D}}\right)$  as
$$\mathcal{F}_{\mathbb{D}}^{\varphi_{\mathbb{D}}}\left(\Omega_{\mathbb{D}}\right)=\mathcal{F}_{\mathbb{D}}^{\varphi_{\mathbb{D}_1}}\left(\Omega_{\mathbb{D}}\right)e_1+\mathcal{F}_{\mathbb{D}}^{\varphi_{\mathbb{D}_2}}\left(\Omega_{\mathbb{D}}\right)e_2,$$
where each of  $\mathcal{F}_{\mathbb{D}}^{\varphi_{\mathbb{D}_1}}\left(\Omega_{\mathbb{D}}\right)$ and $\mathcal{F}_{\mathbb{D}}^{\varphi_{\mathbb{D}_2}}\left(\Omega_{\mathbb{D}}\right)$ are  $\mathbb{C}(i)$ or $\mathbb{C}(j)$ -linear spaces.  
\end{remark}

\begin{proposition}
Let $\left(\Omega_{\mathbb{D}},\Sigma_{\mathbb{D}},\mu_{\mathbb{D}}\right)$ be a $\mathbb{D}$-finite measure space. Then
\begin{itemize}
\item[(i)]$L_{\mathbb{D}}^1(\mu_{\mathbb{D}})=\bigcup\left\{\mathcal{F}_{\mathbb{D}}^{\varphi_{\mathbb{D}}}\left(\Omega_{\mathbb{D}}\right): \varphi_{\mathbb{D}}~~ \mbox{varies over all}~~ \mathbb{D}- \mbox{valued N-functions.}\right\}$.\\
\item[(ii)] $L^{\infty,\mathbb{D}}(\mu_{\mathbb{D}})=\bigcap\left\{\mathcal{F}_{\mathbb{D}}^{\varphi_{\mathbb{D}}}\left(\Omega_{\mathbb{D}}\right): \varphi_{\mathbb{D}}~~ \mbox{varies over all}~~ \mathbb{D}- \mbox{valued N-functions}.\right\}.$
\end{itemize}
\end{proposition}

\begin{proof}
(i) 
Let $f=f_1e_1+f_2e_2 \in \bigcup\left\{\mathcal{F}_{\mathbb{D}}^{\varphi_{\mathbb{D}}}\left(\Omega_{\mathbb{D}}\right): \varphi_{\mathbb{D}}~~ \mbox{varies over all}~~ \mathbb{D}- \mbox{valued N-functions}\right\}.$ This means that $f=f_1e_1+f_2e_2 \in \mathcal{F}_{\mathbb{D}}^{\varphi_{\mathbb{D}}}\left(\Omega_{\mathbb{D}}\right)$ for some $\varphi_{\mathbb{D}}$. So, $f_1 \in \mathcal{F}_{\mathbb{D}}^{\varphi_{\mathbb{D}_1}}\left(\Omega_{\mathbb{D}}\right)$ and $f_2 \in  \mathcal{F}_{\mathbb{D}}^{\varphi_{\mathbb{D}_2}}\left(\Omega_{\mathbb{D}}\right)$, where $\mathcal{F}_{\mathbb{D}}^{\varphi_{\mathbb{D}_1}}\left(\Omega_{\mathbb{D}}\right)=e_1\mathcal{F}_{\mathbb{D}}^{\varphi_{\mathbb{D}}}\left(\Omega_{\mathbb{D}}\right)$ and $\mathcal{F}_{\mathbb{D}}^{\varphi_{\mathbb{D}_2}}\left(\Omega_{\mathbb{D}}\right)=e_2\mathcal{F}_{\mathbb{D}}^{\varphi_{\mathbb{D}}}\left(\Omega_{\mathbb{D}}\right).$ We find that $f_1 \in L_{\mathbb{D}}^{11}(\Omega_{\mathbb{D}})$ and $f_2 \in L_{\mathbb{D}}^{21}(\Omega_{\mathbb{D}})$, where $L_{\mathbb{D}}^{11}(\Omega_{\mathbb{D}})=e_1L_{\mathbb{D}}^{1}(\Omega_{\mathbb{D}})$ and $L_{\mathbb{D}}^{21}(\Omega_{\mathbb{D}})=e_2L_{\mathbb{D}}^{1}(\Omega_{\mathbb{D}})$ are collection of all $\mathbb{C}$-valued integrable functions. Thus $f=f_1e_1+f_2e_2 \in L_{\mathbb{D}}^{1}(\Omega_{\mathbb{D}}).$ Since this is true for every such $f$ in the union, so
$$\bigcup\left\{\mathcal{F}_{\mathbb{D}}^{\varphi_{\mathbb{D}}}\left(\Omega_{\mathbb{D}}\right): \varphi_{\mathbb{D}}~~ \mbox{varies over all}~~ \mathbb{D}- \mbox{valued N-functions.}\right\} \subseteq L_{\mathbb{D}}^1(\Omega_{\mathbb{D}}).$$ Conversely suppose that $f=f_1e_1+f_2e_2 \in L_{\mathbb{D}}^1(\Omega_{\mathbb{D}}).$ Then clearly $f \in \mathcal{A}.$ By using Theorem \ref{th:38}, we can find a $\mathbb{D}$-valued N-function $\varphi_{\mathbb{D}}'=\varphi_{\mathbb{D}_1}'e_1+ \varphi_{\mathbb{D}_2}'e_2$ such that 
$$\int_{\Omega_{\mathbb{D}}}\varphi_{\mathbb{D}}'(|f|_k)d\mu_{\mathbb{D}} <' \infty.$$
Thus $f \in \mathcal{F}_{\mathbb{D}}^{\varphi_{\mathbb{D}}'}\left(\Omega_{\mathbb{D}}\right) \subseteq \bigcup\left\{\mathcal{F}_{\mathbb{D}}^{\varphi_{\mathbb{D}}}\left(\Omega_{\mathbb{D}}\right): \varphi_{\mathbb{D}}~~ \mbox{varies over all}~~ \mathbb{D}- \mbox{valued N-functions.}\right\}.$ This completes the proof.\\
\\
(ii) Let $f=f_1e_1+f_2e_2 \in L^{\infty,\mathbb{D}}(\Omega_{\mathbb{D}})$. Then $f_1 \in L^{\infty,e_1\mathbb{D}}(\Omega_{\mathbb{D}})$ and $f_2 \in L^{\infty,e_2\mathbb{D}}(\Omega_{\mathbb{D}}).$ where $L^{\infty,e_1\mathbb{D}}(\Omega_{\mathbb{D}})=e_1L^{\infty,\mathbb{D}}(\Omega_{\mathbb{D}})$ and $L^{\infty,e_2\mathbb{D}}(\Omega_{\mathbb{D}})=e_2L^{\infty,\mathbb{D}}(\Omega_{\mathbb{D}})$ are the collection of all $\mathbb{C}$-valued essentially bounded functions. Since $L^{\infty,\mathbb{R}} \subseteq \mathcal{F}_{\mathbb{D}}^{\varphi_{\mathbb{D}}}\left(\Omega_{\mathbb{D}}\right)$ for every $\mathbb{R}$-valued N function, we find that $f_1 \in \mathcal{F}_{\mathbb{D}}^{\varphi_{\mathbb{D}_1}}\left(\Omega_{\mathbb{D}}\right)$ and $f_2 \in \mathcal{F}_{\mathbb{D}}^{\varphi_{\mathbb{D}_2}}\left(\Omega_{\mathbb{D}}\right)$  for every $\mathbb{R}$-valued N-functions $\varphi_{\mathbb{D}_1}=e_1 \varphi_{\mathbb{D}}$ and $\varphi_{\mathbb{D}_2}=e_2\varphi_{\mathbb{D}}$,where $\varphi_{\mathbb{D}}$ is a $\mathbb{D}$-valued N-function. So clearly we see that $f=f_1e_1+f_2e_2 \in \mathcal{F}_{\mathbb{D}}^{\varphi_{\mathbb{D}}}\left(\Omega_{\mathbb{D}}\right)$ for every $\mathbb{D}$  -valued N functions. Conversely, suppose that $f=f_1e_1+f_2e_2 \in \bigcap\left\{\mathcal{F}_{\mathbb{D}}^{\varphi_{\mathbb{D}}}\left(\Omega_{\mathbb{D}}\right): \varphi_{\mathbb{D}}~~ \mbox{varies over all}~~ \mathbb{D}- \mbox{valued N-functions}\right\}.$ Then $f=f_1e_1+f_2e_2 \in \mathcal{F}_{\mathbb{D}}^{\varphi_{\mathbb{D}}}\left(\Omega_{\mathbb{D}}\right)$ for every $\mathbb{D}$-valued N functions $\varphi_{\mathbb{D}}$, so that $f_1 \in \mathcal{F}_{\mathbb{D}}^{\varphi_{\mathbb{D}_1}}\left(\Omega_{\mathbb{D}}\right) \subseteq L^{\infty,e_1\mathbb{D}}(\Omega_{\mathbb{D}})$ and  $f_2 \in \mathcal{F}_{\mathbb{D}}^{\varphi_{\mathbb{D}_2}}\left(\Omega_{\mathbb{D}}\right) \subseteq L^{\infty,e_2\mathbb{D}}(\Omega_{\mathbb{D}}).$
Thus $f=f_1e_1+f_2e_2 \in L^{\infty,\mathbb{D}}(\Omega_{\mathbb{D}})=e_1L^{\infty,e_1\mathbb{D}}(\Omega_{\mathbb{D}})+e_2L^{\infty,e_2\mathbb{D}}(\Omega_{\mathbb{D}}).$ This completes the proof.
\end{proof}

\begin{definition}
Let $\mathcal{F}_{\mathbb{D}}^{\varphi_{\mathbb{D}}}\left(\Omega_{\mathbb{D}},\Sigma_{\mathbb{D}},\mu_{\mathbb{D}}\right)$ be the collection as defined above. Then the space $\mathcal{L}_{\mathbb{D}}^{\varphi_{\mathbb{D}}}\left(\Omega_{\mathbb{D}},\Sigma_{\mathbb{D}},\mu_{\mathbb{D}}\right)$ of all $\mathbb{BC}$-valued $\Sigma_{\mathbb{D}}$-measurable function $f: \Omega_{\mathbb{D}} \to \mathbb{BC}$ such that
$$\alpha f \in \mathcal{F}_{\mathbb{D}}^{\varphi_{\mathbb{D}}}\left(\Omega_{\mathbb{D}},\Sigma_{\mathbb{D}},\mu_{\mathbb{D}}\right)$$
$$i.e.,~~~~~\mathcal{L}_{\mathbb{D}}^{\varphi_{\mathbb{D}}}\left(\Omega_{\mathbb{D}},\Sigma_{\mathbb{D}},\mu_{\mathbb{D}}\right)=\left\{f:\Omega_{\mathbb{D}} \to \mathbb{BC}, measurable~: \int_{\Omega_{\mathbb{D}}} \varphi_{\mathbb{D}}(|\alpha f|_k)d\mu_{\mathbb{D}} \leq' \infty \right\}$$ 
for some $\mathbb{D}$-small bicomplex number $\alpha >' 0$ is called the $\mathbb{BC}$-Orlicz space. 
The $\mathbb{BC}$-Orlicz space can further be decompose in two pairs of classical Orlicz spaces and we can write 
$$\mathcal{L}_{\mathbb{D}}^{\varphi_{\mathbb{D}}}\left(\Omega_{\mathbb{D}},\right)= \mathcal{L}_{e_1\mathbb{D}}^{\varphi_{\mathbb{D}_1}}\left(\Omega_{\mathbb{D}},\Sigma_{\mathbb{D}},\mu_{\mathbb{D}}^1\right)e_1+\mathcal{L}_{e_2\mathbb{D}}^{\varphi_{\mathbb{D}_2}}\left(\Omega_{\mathbb{D}},\Sigma_{\mathbb{D}},\mu_{\mathbb{D}}^2\right)$$
\end{definition}

Now we shall define a $\mathbb{D}$-valued modular function denoted by $I_{\varphi_{\mathbb{D}}}^{\mathbb{D}}$ on the hyperbolic Orlicz space as follows
\begin{eqnarray*}
I_{\varphi_{\mathbb{D}}}^{\mathbb{D}}(f)&=& \int_{\Omega_{\mathbb{D}}} \varphi_{\mathbb{D}}(|f|_k)d\mu_{\mathbb{D}} \end{eqnarray*} and by using equation \eqref{eq:23} we  further have

\begin{eqnarray*}
I_{\varphi_{\mathbb{D}}}^{\mathbb{D}}(f)&=&I_{\varphi_{\mathbb{D}_1}}^{\mathbb{D}_1}(e_1f)e_1+I_{\varphi_{\mathbb{D}_2}}^{\mathbb{D}_2}(e_2f)e_2.						
\end{eqnarray*}

\begin{definition}
A $\mathbb{D}$-valued norm on the $\mathbb{BC}$-Orlicz space $\mathcal{L}_{\mathbb{D}}^{\varphi_{\mathbb{D}}}\left(\Omega_{\mathbb{D}},\Sigma_{\mathbb{D}},\mu_{\mathbb{D}}\right)$ is defined as 
$$\|f\|_{\mathbb{D}}^{\varphi_{\mathbb{D}}}=\inf \left\{\alpha~>' 0,: I_{\varphi_{\mathbb{D}}}^{\mathbb{D}}\left(\frac{|f|_k}{\alpha}\right) \leq' 1\right\}$$
Note that such  $\alpha$ will always lie outside the null cone $\mathcal{NC}$ of $\mathbb{D}.$
We say that this is a $\mathbb{D}$-valued Luxemburg norm on $\mathcal{L}_{\mathbb{D}}^{\varphi_{\mathbb{D}}}\left(\Omega_{\mathbb{D}},\Sigma_{\mathbb{D}},\mu_{\mathbb{D}}\right).$ This norm can further be decompose into the following
\begin{eqnarray*}
\|f\|_{\mathbb{D}}^{\varphi_{\mathbb{D}}} & = & \|e_1f\|_{\mathbb{D}_1}^{\varphi_{\mathbb{D}_1}}e_1+\|e_2f\|_{\mathbb{D}_2}^{\varphi_{\mathbb{D}_2}}e_1,
\end{eqnarray*}
where $\|.\|_{\mathbb{D}_1}^{\varphi_{\mathbb{D}_1}}$ and $\|.\|_{\mathbb{D}_2}^{\varphi_{\mathbb{D}_2}}$ are the real valued Luxemburg norm on the classical Orlicz spaces $\mathcal{L}_{e_1\mathbb{D}}^{\varphi_{\mathbb{D}_1}}\left(\Omega_{\mathbb{D}},\Sigma_{\mathbb{D}},\mu_{\mathbb{D}}^1\right)$ and $\mathcal{L}_{e_2\mathbb{D}}^{\varphi_{\mathbb{D}_2}}\left(\Omega_{\mathbb{D}},\Sigma_{\mathbb{D}},\mu_{\mathbb{D}}^2\right)$ defined by using the $\mathbb{R}$-valued convex functions $\varphi_{\mathbb{D}_1}$  and $\varphi_{\mathbb{D}_2}$ respectively.
\end{definition}

\begin{theorem}
The $\mathbb{BC}$-linear module $\mathcal{L}_{\mathbb{D}}^{\varphi_{\mathbb{D}}}\left(\Omega_{\mathbb{D}},\Sigma_{\mathbb{D}},\mu_{\mathbb{D}}\right)$ with the $\mathbb{D}$-Luxemburg norm $\|.\|_{\mathbb{D}}^{\varphi_{\mathbb{D}}}$ is a complete $\mathbb{D}$-norm vector module.
\end{theorem}

\begin{proof}
Proof of this theorem can be done in a similar way as was done for the completeness of the $\mathbb{D}$-metric spaces. 
\end{proof}

\begin{theorem}
Let $\mathcal{L}_{\mathbb{D}}^{\varphi_{\mathbb{D}}}\left(\Omega_{\mathbb{D}},\Sigma_{\mathbb{D}},\mu_{\mathbb{D}}\right)$ be a $\mathbb{BC}$-Orlicz space with the $\mathbb{D}$-Luxemburg norm $\|.\|_{\mathbb{D}}^{\varphi_{\mathbb{D}}}$. Then prove that 
$$\|f\|_{\mathbb{D}}^{\varphi_{\mathbb{D}}} \leq' 1~~~~\mbox{if and only if}~~~~I_{\varphi_{\mathbb{D}}}^{\mathbb{D}}(f) \leq' 1. $$ 
\end{theorem}
\begin{proof}
First suppose that 
$\|f\|_{\mathbb{D}}^{\varphi_{\mathbb{D}}} \leq' 1,$ Thus 
\begin{eqnarray}
\label{eq:30}
\|e_1f\|_{\mathbb{D}_1}^{\varphi_{\mathbb{D}_1}}e_1+\|e_2f\|_{\mathbb{D}_2}^{\varphi_{\mathbb{D}_2}}e_2 \leq' 1.
\end{eqnarray}
 Multiplying equation $\eqref{eq:30}$ by $e_1$ and $e_2$ we get $\|e_1f\|_{\mathbb{D}_1}^{\varphi_{\mathbb{D}_1}}e_1 \leq' e_1$ and $\|e_1f\|_{\mathbb{D}_1}^{\varphi_{\mathbb{D}_1}}e_1 \leq' e_1.$ By the partial ordering relation we get
$\|e_1f\|_{\mathbb{D}_1}^{\varphi_{\mathbb{D}_1}} \leq 1$ and  $\|e_2f\|_{\mathbb{D}_2}^{\varphi_{\mathbb{D}_2}} \leq 1.$ Since the result is true for $\mathbb{R}$-Orlicz space, we have
\begin{equation}
\label{eq:31}
I_{\varphi_{\mathbb{D}_1}}^{\mathbb{D}_1}(e_1f) \leq 1~\mbox{and}~ I_{\varphi_{\mathbb{D}_2}}^{\mathbb{D}_2}(e_2f) \leq 1. 
\end{equation}
From equation \eqref{eq:31} we get 
$I_{\varphi_{\mathbb{D}_1}}^{\mathbb{D}_1}(e_1f)e_1 + I_{\varphi_{\mathbb{D}_2}}^{\mathbb{D}_2}(e_2f)e_2 \leq' 1.$ This we get $I_{\varphi_{\mathbb{D}}}^{\mathbb{D}}(f) \leq' 1.$ Similarly converse can be retraced back.
\end{proof}

\begin{theorem}
Let $\left(\varphi_{\mathbb{D}} , \psi_{\mathbb{D}}\right)$ be a $\mathbb{D}$-complementary pair of Young functions and let $\mathcal{L}_{\mathbb{D}}^{\varphi_{\mathbb{D}}}\left(\Omega_{\mathbb{D}},\Sigma_{\mathbb{D}},\mu_{\mathbb{D}}\right)$ and $\mathcal{L}_{\mathbb{D}}^{\psi_{\mathbb{D}}}\left(\Omega_{\mathbb{D}},\Sigma_{\mathbb{D}},\mu_{\mathbb{D}}\right)$ be the corresponding $\mathbb{BC}$-Orlicz spaces with $\mathbb{D}$-valued Luxemburg norm $\|.\|_{\mathbb{D}}^{\varphi_{\mathbb{D}}}$ and $\|.\|_{\mathbb{D}}^{\psi_{\mathbb{D}}}$ respectively. If $f \in \mathcal{L}_{\mathbb{D}}^{\varphi_{\mathbb{D}}}\left(\Omega_{\mathbb{D}},\Sigma_{\mathbb{D}},\mu_{\mathbb{D}}\right)$ and $g \in \mathcal{L}_{\mathbb{D}}^{\psi_{\mathbb{D}}}\left(\Omega_{\mathbb{D}},\Sigma_{\mathbb{D}},\mu_{\mathbb{D}}\right),$ then one has 
\begin{equation}
\label{eq:Holder1}
\int_{\Omega_{\mathbb{D}}}|fg|_kd\mu_{\mathbb{D}} \leq ' 2\|f\|_{\mathbb{D}}^{\varphi_{\mathbb{D}}} \cdot \|g\|_{\mathbb{D}}^{\psi_{\mathbb{D}}}.
\end{equation}
\end{theorem}

\begin{proof}
If either $\|f\|_{\mathbb{D}}^{\varphi_{\mathbb{D}}}$ or $\|g\|_{\mathbb{D}}^{\psi_{\mathbb{D}}}$ is zero so that $f=0$ or $g=0$ $a.e$, then it is clear that \eqref{eq:Holder1} is trivially true. So let us assume that $\|f\|_{\mathbb{D}}^{\varphi_{\mathbb{D}}}  >' 0$ and $\|g\|_{\mathbb{D}}^{\psi_{\mathbb{D}}} >' 0.$ Then by Young's inequality, we have 
\begin{equation}
\label{eq:Holder2}
\frac{|fg|_k}{\|f\|_{\mathbb{D}}^{\varphi_{\mathbb{D}}} \cdot \|g\|_{\mathbb{D}}^{\psi_{\mathbb{D}}}} \leq ' \varphi_{\mathbb{D}} \left( \frac{|f|_k}{\|f\|_{\mathbb{D}}^{\varphi_{\mathbb{D}}}}\right) + \psi_{\mathbb{D}}\left( \frac{|g|_k}{\|g\|_{\mathbb{D}}^{\psi_{\mathbb{D}}}}\right)
\end{equation} 
Integrating equation \eqref{eq:Holder2}, we get
\begin{equation}
\label{eq:Holder3}
\int_{\Omega_{\mathbb{D}}}\frac{|fg|_k}{\|f\|_{\mathbb{D}}^{\varphi_{\mathbb{D}}} \cdot \|g\|_{\mathbb{D}}^{\psi_{\mathbb{D}}}} d\mu_{\mathbb{D}} \leq ' \int_{\Omega_{\mathbb{D}}} \varphi_{\mathbb{D}} \left( \frac{|f|_k}{\|f\|_{\mathbb{D}}^{\varphi_{\mathbb{D}}}}\right)d\mu_{\mathbb{D}} + \int_{\Omega_{\mathbb{D}}} \psi_{\mathbb{D}}\left( \frac{|g|_k}{\|g\|_{\mathbb{D}}^{\psi_{\mathbb{D}}}}\right)d\mu_{\mathbb{D}}.
\end{equation} 
Now since we know that the first and second terms on the right hand side of equation \eqref{eq:Holder3} is less than or equal to one. Thus we 
find that 

\begin{eqnarray*}
\int_{\Omega_{\mathbb{D}}}|fg|_k d\mu_{\mathbb{D}} & \leq ' & 2 \|f\|_{\mathbb{D}}^{\varphi_{\mathbb{D}}} \cdot \|g\|_{\mathbb{D}}^{\psi_{\mathbb{D}}}.
\end{eqnarray*}

\end{proof}

\section{\bf Bicomplex linear operators and spectrum of bicomplex unilateral shift operators with $\mathbb{D}$-norm}

The $\mathbb{BC}$-linear operators are discussed in detail in  \cite{Al-Lu-Sh-Sr13}.  The spectrum, resolvent and resolvent operators were   defined for a bicomplex linear operators in \cite{Co-Sa-St13}.  Here in this section we show that  point spectrum of the unilateral shift operators on $l^2(\mathbb{BC})$ is the null cone. 

\begin{definition}
Let $\left(\Omega_{\mathbb{D}},\Sigma_{\mathbb{D}},\mu_{\mathbb{D}}\right)$ be a $\sigma$-finite $\mathbb{D}$-measure space and let $T_{\mathbb{D}}:\Omega_{\mathbb{D}} \to \Omega_{\mathbb{D}}$ be a $\Sigma_{\mathbb{D}}$-measurable transformation. $i.e.,~~T_{\mathbb{D}}^{-1}(E) \in \Sigma_{\mathbb{D}}$ for any $E \in \Sigma_{\mathbb{D}}$ then $T_{\mathbb{D}}$ is said to be a non-singular $\Sigma_{\mathbb{D}}$-measurable transformation if $\mu_{\mathbb{D}}(T_{\mathbb{D}}^{-1}(E))=0$ for every $E \in \Sigma_{\mathbb{D}}$ with $\mu_{\mathbb{D}}(E)=0.$
A non-singular measurable transformation $T_{\mathbb{D}}$ induces a $\mathbb{BC}$-linear operator $C_T$ from $\mathcal{M}^o(\left(\Omega_{\mathbb{D}},\Sigma_{\mathbb{D}},\mu_{\mathbb{D}}\right))$ into itself defined by
$$C_{T_{\mathbb{D}}}f=foT_{\mathbb{D}}~~~~\forall~~~f \in \mathcal{M}^o(\left(\Omega_{\mathbb{D}},\Sigma_{\mathbb{D}},\mu_{\mathbb{D}}\right)).$$
where $\mathcal{M}^o(\left(\Omega_{\mathbb{D}},\Sigma_{\mathbb{D}},\mu_{\mathbb{D}}\right))$, denotes the collection of all equivalence class of $\mathbb{BC}$-valued $\Sigma_{\mathbb{D}}$-measurable functions which are identified $\mu_{\mathbb{D}}$-almost everywhere.

\end{definition}
\begin{remark} The non-singularity of $T_{\mathbb{D}}$ ensures that $C_{T_{\mathbb{D}}}$ is well defined as a $\mathbb{BC}$-linear operator from $\mathcal{M}^o(\left(\Omega_{\mathbb{D}},\Sigma_{\mathbb{D}},\mu_{\mathbb{D}}\right))$ into $\mathcal{M}^o(\left(\Omega_{\mathbb{D}},\Sigma_{\mathbb{D}},\mu_{\mathbb{D}}\right)).$
\end{remark}

Now since $f\in \mathcal{M}^o(\left(\Omega_{\mathbb{D}},\Sigma_{\mathbb{D}},\mu_{\mathbb{D}}\right))$ implies that $f=f_1e_1+f_2e_2$ and so we can decompose the composition operator $C_{T_{\mathbb{D}}}$ induced by  $T_{\mathbb{D}}$ as 
\begin{eqnarray*}
C_{T_{\mathbb{D}}}(f(t)) & = & C_{T_{\mathbb{D}}}(f_1(t)e_1+f_2(t)e_2)\\
                         & = & C_{T_{\mathbb{D}}}(f_1(t))e_1+C_{T_{\mathbb{D}}}(f_2(t))e_2\\												
i.e.,  C_{T_{\mathbb{D}}}(f) & = & C_{T_{\mathbb{D}}}(f_1)e_1+C_{T_{\mathbb{D}}}(f_2)e_2\\ 
\end{eqnarray*}
\begin{remark}
$C_{T_{\mathbb{D}}}$ is a composition operator on $\mathcal{M}^o(\left(\Omega_{\mathbb{D}},\Sigma_{\mathbb{D}},\mu_{\mathbb{D}}\right))$ if and only if $C_{T_{\mathbb{D}}}$ is a composition operator on $\mathcal{M}^o(\left(\Omega_{\mathbb{D}},\Sigma_{\mathbb{D}},\mu_{\mathbb{D}}^1\right))$ and $\mathcal{M}^o(\left(\Omega_{\mathbb{D}},\Sigma_{\mathbb{D}},\mu_{\mathbb{D}}^2\right)).$
\end{remark}
Next we shall define multiplication operator on $\mathcal{M}^o(\left(\Omega_{\mathbb{D}},\Sigma_{\mathbb{D}},\mu_{\mathbb{D}}\right)).$

\begin{definition}
Let $\psi_{\mathbb{D}}:\Omega_{\mathbb{D}} \to \mathbb{BC}$ be a $\mathbb{BC}$-valued $\Sigma_{\mathbb{D}}$- measurable function. We define a multiplication operator 
$M_{\psi_{\mathbb{D}}}:\mathcal{M}^o(\left(\Omega_{\mathbb{D}},\Sigma_{\mathbb{D}},\mu_{\mathbb{D}}\right)) \to \mathcal{M}^o(\left(\Omega_{\mathbb{D}},\Sigma_{\mathbb{D}},\mu_{\mathbb{D}}\right))$ by 
$$M_{\psi_{\mathbb{D}}}(f(t))=(\psi_{\mathbb{D}} \cdot f)(t)=\psi_{\mathbb{D}}(t)f(t)$$
for all $t \in \Omega_{\mathbb{D}}$ and $f=f_1e_1+f_2e_2 \in \mathcal{M}^o(\left(\Omega_{\mathbb{D}},\Sigma_{\mathbb{D}},\mu_{\mathbb{D}}\right)).$
\end{definition}

\begin{remark}
It is easy to see that as $\psi_{\mathbb{D}},f$ are $\mathbb{BC}$-valued $\Sigma_{\mathbb{D}}$-measurable functions implies that $\psi_{\mathbb{D}}(t)=\psi_{\mathbb{D}}^1(t)e_1+\psi_{\mathbb{D}}^2(t)e_2$ and $f(t)=f_1(t)e_1+f_2(t)e_2$ so that 
\begin{eqnarray*}
M_{\psi_{\mathbb{D}}}(f(t)) & = & (\psi_{\mathbb{D}} \cdot f)(t)\\
                            & = & \psi_{\mathbb{D}}(t) \cdot f(t)\\
														& = & (\psi_{\mathbb{D}}^1(t)e_1+\psi_{\mathbb{D}}^2(t)e_2) \cdot (f_1(t)e_1+f_2(t)e_2)\\
														& = & (\psi_{\mathbb{D}}^1(t)\cdot  f_1(t))e_1+(\psi_{\mathbb{D}}^2(t)\cdot f_2(t))e_2\\
														& = & M_{\psi_{\mathbb{D}}^1}(f_1(t))e_1+M_{\psi_{\mathbb{D}}^2}(f_2(t))e_2.\\
	i.e.,M_{\psi_{\mathbb{D}}}(f(t)) & = & M_{\psi_{\mathbb{D}}^1}(f_1(t))e_1+M_{\psi_{\mathbb{D}}^2}(f_2(t))e_2.\\
	\end{eqnarray*}
Thus we see that $M_{\psi_{\mathbb{D}}}(f(t))$ is a multiplication operator induced by $\psi_{\mathbb{D}}$ on $\mathcal{M}^o(\left(\Omega_{\mathbb{D}},\Sigma_{\mathbb{D}},\mu_{\mathbb{D}}\right))$ if and only if $M_{\psi_{\mathbb{D}}^1}$ and $M_{\psi_{\mathbb{D}}^2}$ are multiplication operators induced by $\psi_{\mathbb{D}}^1$ and $\psi_{\mathbb{D}}^2$ on $\mathcal{M}^o(\left(\Omega_{\mathbb{D}},\Sigma_{\mathbb{D}},\mu_{\mathbb{D}}^1\right))$ and $\mathcal{M}^o(\left(\Omega_{\mathbb{D}},\Sigma_{\mathbb{D}},\mu_{\mathbb{D}}^2\right))$ respectively.      
\end{remark}

\begin{definition}
\label{def:Spec1}
Let $X$ be a $\mathbb{D}$-normed $\mathbb{BC}$-module and $T_{\mathbb{D}}:D(T_{\mathbb{D}})\subset X \to X$ be a $\mathbb{BC}$-linear operator. Then we define the resolvent set $\rho_{\mathbb{D}}(T_{\mathbb{D}})$ of $T_{\mathbb{D}}$ to be the set of all $\lambda \in \mathbb{BC}$ for which the following conditions are satisfied.
\begin{itemize}
\item[(i)] $(T_{\mathbb{D}} -\lambda I)^{-1}$ exists.\\
\item[(ii)] $(T_{\mathbb{D}} -\lambda I)^{-1}$ is $\mathbb{D}$-bounded.\\
\item[(iii)] $(T_{\mathbb{D}} -\lambda I)^{-1}$ is defined on a set which is dense in $X$.\\
\end{itemize}
\end{definition}

\begin{definition}
The spectrum of $T_{\mathbb{D}}$ denoted by $\sigma_{\mathbb{D}}(T_{\mathbb{D}})$ is defined as the complement of $\rho_{\mathbb{D}}(T_{\mathbb{D}})$.
$$i.e.,~~~~\rho_{\mathbb{D}}(T_{\mathbb{D}})=\left\{\lambda \in \mathbb{BC}: (T_{\mathbb{D}} -\lambda I)~~~\mbox{is not invertible}\right\}.$$ 
Further we define an operator $R(\lambda :T_{\mathbb{D}})=(T_{\mathbb{D}} -\lambda I)^{-1}$ for every $\lambda \in \rho_{\mathbb{D}}(T_{\mathbb{D}})$  called the resolvent operator.
\end{definition}

Now the spectrum of a $\mathbb{D}$-bounded $\mathbb{BC}$-linear operator can further be partitioned into three disjoint sets as follows.
\begin{definition}
The point spectrum or the discrete spectrum denoted by $\sigma_{\mathbb{D} p} (T_{\mathbb{D}})$ is the set of all bicomplex numbers  $\lambda$ such that $(T_{\mathbb{D}} -\lambda I)^{-1}$ does not exist.
\end{definition}

\begin{definition}
The residual spectrum denoted by $\sigma_{\mathbb{D}r} (T_{\mathbb{D}})$ is the set of all bicomplex number $\lambda$ such that 
$(T_{\mathbb{D}} -\lambda I)^{-1}$ exists (may or may not be $\mathbb{D}$-bounded), but domain of $(T_{\mathbb{D}} -\lambda I)^{-1}$ is not dense in $X$.
\end{definition}

\begin{definition}
The continuous spectrum denoted by $\sigma_{\mathbb{D}c} (T_{\mathbb{D}})$ is the set of all bicomplex number $\lambda$ such that
$(T_{\mathbb{D}} -\lambda I)^{-1}$ exists and is dense in $X$ but $(T_{\mathbb{D}} -\lambda I)^{-1}$ is not $\mathbb{D}$-bounded.
\end{definition}
Thus we can write $\mathbb{D}= \rho_{\mathbb{D}}(T_{\mathbb{D}})\cup \sigma_{\mathbb{D}}(T_{\mathbb{D}})=\rho_{\mathbb{D}}(T_{\mathbb{D}})\cup\sigma_{\mathbb{D} p} (T) \cup \sigma_{\mathbb{D}r} (T) \cup \sigma_{\mathbb{D}c}(T).$

\begin{definition}
Let $X$ and $Y$ be a $\mathbb{D}$-normed $\mathbb{BC}$-module and $T_{\mathbb{D}}:X \to Y$ be a $\mathbb{BC}$-linear operator. Then the adjoint $T^*:X^* \to Y^*$ defined as 
$$T^*(f)(x)=f(T(x))~~~~\forall~~~x~~ \in X ~\mbox{and}~ f \in Y^*.$$
\end{definition}
\begin{lemma}
Let $l^p(\mathbb{BC})$ be the space all $p$-summable sequence $x=(x_1,x_2,x_3,....)$ of bicomplex numbers such that 
\begin{equation}
\sum_{i=1}^{\infty}|x_i|_k^p =|x_1|_k^p+|x_2|_k^p+|x_3|_k^p+ \cdots~\mbox{converges in}~\mathbb{D}.
\end{equation} Then $\sum_{i=1}^{\infty}|x_i|_k^p <' \infty$ if and only if $\sum_{i=1}^{\infty} |\xi_1^i|^p < \infty$ and $\sum_{i=1}^{\infty} |\xi_2^i|^p < \infty$ where $x_i=\xi_1^i e_1+\xi_2^i e_2.$ 
\end{lemma}
The $\mathbb{D}$-norm on $l^p(\mathbb{BC})$ is defined as $$\|x\|_{\mathbb{D}}=\left(\sum_{i=1}^{\infty}|x_i|_k^p\right)^{\frac{1}{p}}.$$

\begin{proposition}
\label{Pr:1}
Let $l^2(\mathbb{BC})$ be the space of all square summable sequence of bicomplex numbers. Then an element $x=(x_i)$ in $l^2(\mathbb{BC})$ is in $\mathcal{NC}$ of $l^2(\mathbb{BC})$ if and only if $x_i$ is either zero or $x_i \in \mathcal{NC}$ of $\mathbb{BC}$ for some $i,$ where $\mathcal{NC}$ denotes the null cone.
\end{proposition}

\begin{proof}
First suppose that $x=(x_i)$ belongs to the null cone of $l^2(\mathbb{BC}).$ Then there exists $0 \ne y=(y_i)$ in $l^2(\mathbb{BC})$ such that $xy=0.$ Thus one has 
$$(x_1y_1,x_2y_2,x_3y_3, \cdots)=(0,0,0,\cdots).$$
Now since $y \ne 0$ implies that $y_i \ne 0$ for some $i$ say $y_k \ne 0$. Therefore $x_ky_k=0$ and $y_k \ne 0$ implies that either $x_k=0$ or  $x_k$ has to be an element of the null cone $\mathcal{NC}$ of $\mathbb{BC}.$ This proves the direct part.\\
Conversly suppose $x=(x_i) \in l^2(\mathbb{BC})$ be such that $x_i=0$  or $x_i \in \mathcal{NC}$ for some i, implies that $x_iy_i=0$ for some $y_i \in \mathbb{BC}$. and so that $y=(0,0,0 \cdots ,y_i, \cdots 0,0,0) \in l^2(\mathbb{BC})$ such that $xy=0$ implies that $x$ is in the null cone.
\end{proof}

In \cite{Co-Sa-St13}, it was proved that the spectrum of any bicomplex linear operator is unbounded. Here we show that the spectrum of unilateral shift operators is the null cone.  

\begin{example}
Define the unilateral forward shift operator $S_{\mathbb{D}}:l^2(\mathbb{BC}) \to l^2(\mathbb{BC})$ by
$$S_{\mathbb{D}}(x_1,x_2,x_3, \cdots)=(0,x_1,x_2,x_3, \cdots)~~~~\forall~~{x_i} \in l^2(\mathbb{BC}).$$
Let $(x_1,x_2,x_2, \cdots) \in l^2(\mathbb{BC})$ then we see that
\begin{eqnarray*} 
\|S_{\mathbb{D}}x\|_{\mathbb{D}}^2 & = &\|S_{\mathbb{D}}(x_1,x_2,x_3, \cdots)\|_{\mathbb{D}}^2\\
                                  & = & \|(0,x_1,x_2,x_3, \cdots)\|_{\mathbb{D}}^2\\
																	& = & \sum_{n=1}^{\infty}|x_i|_k^2\\
																	& = & \|x\|_{\mathbb{D}}^2.
\end{eqnarray*}
Thus we find that $\|S_{\mathbb{D}}\|_{\mathbb{D}}=1$ which shows that $S_{\mathbb{BC}}$ is a bounded linear operator.\\
Further for $\lambda=0$   the inverse map$(S_{\mathbb{D}}-\lambda I)^{-1}=S_{\mathbb{D}}^{-1}:S_{\mathbb{D}}(l^2(\mathbb{D}))$ exists and is the unilateral left shift operator defined by 
$$S_{\mathbb{D}}^{-1}(x_1,x_2,x_3, \cdots)=(x_2,x_3,x_4, \cdots).$$
But $S_{\mathbb{D}}^{-1}$ is not defined on a dense subset of $l^2(\mathbb{D})$ as $S_{\mathbb{D}}(l^2(\mathbb{D}))$ is clearly not dense in $l^2(\mathbb{BC})$ on which $S_{\mathbb{D}}^{-1}$ is defined. Thus $\lambda=0 \in  \sigma_{\mathbb{D}r}(S_{\mathbb{D}})$. Further for $\lambda \ne 0 \in \sigma_{\mathbb{D} p} (S_{\mathbb{D}})$ we have a non zero $x=(x_1,x_2,x_3 \cdots) \in l^2(\mathbb{BC})$ such that
\begin{eqnarray*}
(S_{\mathbb{D}}-\lambda I)x & = & 0\\
\Leftrightarrow (S_{\mathbb{D}}-\lambda I)(x_1,x_2,x_3 \cdots) & = & (0,0,0 \cdots)\\
\Leftrightarrow S_{\mathbb{D}}(x_1,x_2,x_3 \cdots)-\lambda I(x_1,x_2,x_3 \cdots) & = & (0,0,0 \cdots)\\
\Leftrightarrow (0,x_1,x_2, \cdots)-\lambda (x_1,x_2,x_3 \cdots) & = & (0,0,0 \cdots)\\
\Leftrightarrow (-\lambda x_1,x_1-\lambda x_2,x_2-\lambda x_3, \cdots) & = & (0,0,0 \cdots)\\
\end{eqnarray*} 
\begin{equation}
\label{eqaa:0}
   \left(-\lambda x_1,x_1-\lambda x_2,x_2-\lambda x_3, \cdots\right)  =  (0,0,0 \cdots)
\end{equation}
Now if $x=(x_1,x_2,x_3 \cdots) \notin \mathcal{NC}$, then we see that $x_i \ne 0$ and $x_i \notin \mathcal{NC}$ of $\mathbb{BC}$ for all $i.$ Then by equation \eqref{eqaa:0} we have $\lambda x_1=0$ implies that $\lambda \in \mathcal{NC}.$ 
Thus we see that $0 \ne \lambda \in \sigma_{\mathbb{D} p} (S_{\mathbb{D}}),$ implies that $\lambda \in \mathcal{NC}$. Therefore 
$$\sigma_{\mathbb{D} p} (S_{\mathbb{D}})= \mathcal{NC}.$$  
\end{example}

\begin{example}
Define $S_{\mathbb{D}}:l^2(\mathbb{BC}) \to l^2(\mathbb{BC})$ by
$$S_{\mathbb{D}}(x_1,x_2,x_3, \cdots)=(0,x_1,x_2,x_3, \cdots).$$ In fact $S_{\mathbb{D}}$ is the unilateral right shift operator. Further note that 
$$\left(l^2(\mathbb{BC})\right)^*=l^2(\mathbb{BC}).$$
Now consider $S_{\mathbb{D}}^*:l^2(\mathbb{BC}) \to l^2(\mathbb{BC}).$ Then for $y=f=(y_1,y_2,y_3, \cdots) \in l^2(\mathbb{BC})$ and $x=(x_1,x_2,x_3, \cdots) \in l^2(\mathbb{BC}),$ we have 
\begin{eqnarray*}
(S_{\mathbb{D}}^*f)(x) & = & f(S_{\mathbb{D}}(x))\\
                       & = & f((0,x_1,x_2,x_3, \cdots))\\
											 & = & y_2x_1+y_3x_2+y_4x_3+ \cdots\\
											 & = & \sum_{n=1}^{\infty} y_{n+1}x_n.\\
i.e.,~~~~~(S_{\mathbb{D}}^*f)(x) & = & \sum_{n=1}^{\infty} y_{n+1}x_n.
\end{eqnarray*} 
Therefore we have
\begin{equation}
i.e.,~~~~~S_{\mathbb{D}}^*f = (y_2,y_3,y_4, \cdots)\\
\end{equation}

\begin{equation}
i.e.,~~~~~S_{\mathbb{D}}^*(y_1,y_2,y_3, \cdots)  = (y_2,y_3,y_4, \cdots).
\end{equation}
which is nothing but the unilateral right shift operator.\\
Let us $x=(x_1,x_2,x_3, \cdots) \in l^2(\mathbb{D})$. Then
\begin{eqnarray*}
\|S_{\mathbb{D}}^*(x)\|_{\mathbb{D}}^2 & = &\|(x_2,x_3,x_4, \cdots)\|_{\mathbb{D}}\\
                          & = & \sum_{n=2}^{\infty} |x_n|_k^2\\
													& \leq' & \sum_{n=1}^{\infty} |x_n|_k^2\\
													& = & \|x\|_{\mathbb{D}}^2.
\end{eqnarray*}
Thus we find that $\|S_{\mathbb{D}}^*\|_{\mathbb{D}} \leq'1$ and so $S_{\mathbb{D}}^*$ is bounded linear operator.
Further for $\lambda=0$ 
$$(S_{\mathbb{D}}^*-\lambda I)^{-1}={S_{\mathbb{D}}^*}^{-1}: S_{\mathbb{D}}(l^2(\mathbb{BC})) \to l^2(\mathbb{BC})$$
exists and is the unilateral right shift operator defined by 
$${S_{\mathbb{D}}^*}^{-1}(\xi_1,\xi_2,\xi_3)==(0,\xi_1,\xi_2,\xi_3)$$
and ${S_{\mathbb{D}}^*}$ satisfies both $(i)$ and $(ii)$ of Definition \ref{def:Spec1} and so therefore $0 \in \rho_{\mathbb{D}}(S_{\mathbb{D}}^*).$
\end{example}
\section{ \bf Musielak Orlicz Modules}
Now we shall finally give a short introduction of Musielak Orlicz functions whose range is the set of hyperbolic numbers. Let $\mathbb{D}$, $\mathbb{D}^+$ and $\mathbb{N}$ stands for the set of hyperbolic numbers,non-negative hyperbolic number and the set of natural numbers respectively. Let $\left(\Omega_{\mathbb{D}},\Sigma_{\mathbb{D}},\mu_{\mathbb{D}}\right)$ be a $\sigma$-finite complete $\mathbb{D}$-measure space.

\begin{definition}
A map $\varphi_{\mathbb{D}}:\Omega_{\mathbb{D}} \times \mathbb{D}^+ \to \mathbb{D}^+$ is said to be a $\mathbb{D}$-Musielak Orlicz function if it satisfies the following conditions:
\begin{itemize}
\item[(i)] $\varphi_{\mathbb{D}}(.,u)$ is a $\mathbb{D}$-valued $\Sigma_{\mathbb{D}}$-measurable function for each 
$u \in \mathbb{D}$.\\
\item[(ii)] $\varphi_{\mathbb{D}}(t,.)$ is a $\mathbb{D}$-valued convex function.\\
\item[(iii)] $\varphi_{\mathbb{D}}(.,u)=0 \Leftrightarrow u=0.$\\
\item[(iv)] $\lim_{u \to +\infty}\varphi_{\mathbb{D}}(.,u) = +\infty$.\\
\end{itemize}
where $+\infty=ae_1+\infty e_2=+\infty e_1+b e_2=+\infty e_1++\infty e_1.$
\end{definition}
Let $\mathcal{M}^o(\left(\Omega_{\mathbb{D}},\Sigma_{\mathbb{D}},\mu_{\mathbb{D}}\right))$ be the equivalence classes of all $\mathbb{BC}$-valued $\Sigma_{\mathbb{D}}$-measurable functions which are identified $\mu_{\mathbb{D}}$ almost everywhere and $\varphi_{\mathbb{D}}$  be any $\mathbb{D}$-Musielak Orlicz function. Then we define a $\mathbb{D}$-convex modular function $I_{\varphi_{\mathbb{D}}}^{\mathbb{D}}$ on $\mathcal{M}^o(\left(\Omega_{\mathbb{D}},\Sigma_{\mathbb{D}},\mu_{\mathbb{D}}\right))$ by 
$$I_{\varphi_{\mathbb{D}}}^{\mathbb{D}}(f)=\int_{\Omega_{\mathbb{D}}} \varphi_{\mathbb{D}}(t,|f(t)|_k)d\mu_{\mathbb{D}}.$$
Here $\varphi_{\mathbb{D}}$ is a $\mathbb{D}$-valued convex function. We also assume the homogenity property of $\varphi_{\mathbb{D}}$ and so, we have 
$$\varphi_{\mathbb{D}}(t,|f(t)|_k)= \varphi_{\mathbb{D}_1}(t_1,|f_1(t_1)|)e_1+\varphi_{\mathbb{D}_2}(t_2,|f_2(t_2)|)e_2.$$
Thus
\begin{eqnarray*}
I_{\varphi_{\mathbb{D}}}^{\mathbb{D}}(f) & = & \int_{\Omega_{\mathbb{D}}} \varphi_{\mathbb{D}}(t,|f(t)|_k)d\mu_{\mathbb{D}}\\
                                         & = & \int_{\Omega_{\mathbb{D}}}\left( \varphi_{\mathbb{D}_1}(t_1,|f_1(t_1)|)e_1+\varphi_{\mathbb{D}_2}(t_1,|f_1(t_1)|)e_1\right) d(\mu_{\mathbb{D}}^1 e_1+ \mu_{\mathbb{D}}^2 e_2)\\
& = & \left(\int_{\Omega_{\mathbb{D}}}\varphi_{\mathbb{D}_1}(t_1,|f_1(t_1)|)d\mu_{\mathbb{D}}^1\right) e_1 -\left(\int_{\Omega_{\mathbb{D}}}\varphi_{\mathbb{D}_2}(t_2,|f_2(t_2)|)d\mu_{\mathbb{D}}^2\right)e_2\\
& = & I_{\varphi_{\mathbb{D}_1}}^{\mathbb{D}}(f_1) e_1 -I_{\varphi_{\mathbb{D}_2}}^{\mathbb{D}}(f_2)e_2.\\
\end{eqnarray*}
Here we have 
$$I_{\varphi_{\mathbb{D}}}^{\mathbb{D}}(f) =  I_{\varphi_{\mathbb{D}_1}}^{\mathbb{D}}(f_1) e_1 - I_{\varphi_{\mathbb{D}_2}}^{\mathbb{D}}(f_2)e_2.$$

\begin{definition}
The Musielak Orlicz space generated by the Musielak Orlicz function $\varphi_{\mathbb{D}}(.,.): \Omega_{\mathbb{D}} \times \mathbb{D}^+ \to \mathbb{D}^+$ is the defined as the collection
$$\mathcal{L}_{\mathbb{D}}^{\varphi_{\mathbb{D}}}\left(\Omega_{\mathbb{D}},\Sigma_{\mathbb{D}},\mu_{\mathbb{D}}\right)=\left\{f \in \mathcal{M}^o(\left(\Omega_{\mathbb{D}},\Sigma_{\mathbb{D}},\mu_{\mathbb{D}}\right)) :I_{\varphi_{\mathbb{D}}}^{\mathbb{D}}(\alpha f) <' \infty \right\}$$
where  $I_{\varphi_{\mathbb{D}}}^{\mathbb{D}}$ is the modular function and $\alpha$ is some positive hyperbolic number.
This space is equipped with the $\mathbb{D}$-Luxemburg norm defined as
$$\|f\|_{\varphi_{\mathbb{D}}}^{\mathbb{D}}= \inf\left\{\lambda >' 0: I_{\varphi_{\mathbb{D}}}^{\mathbb{D}}(\frac{f}{\alpha }) \leq' 1 \right\}$$
and such $\lambda$ will always lie in the complement of the null cone of $\mathbb{D}.$
\end{definition}

\bibliographystyle{plain}

\noindent \textit{Department of Mathematics,\; University of Jammu, \;Jammu,  J\&K - 180 006, India.}\\
E-mail :\textit{ romeshmath@gmail.com.com}\\
E-mail :\textit{ kailash.maths@gmail.com}\\
E-mail :\textit{ joneytun123@gmail.com.com}\\
E-mail.:\textit{ shgunnjamwal04@gmail.com}\\
\end{document}